\documentclass[article,10pt]{IEEEtran}
\usepackage[named]{algo}
\usepackage{graphics}
\usepackage{graphicx}
\usepackage{stfloats}
\usepackage{amssymb}
\usepackage{amsmath}
\usepackage{amsfonts}
\usepackage{subfigure}
\usepackage{float}
\usepackage{url}

\newcommand{\Section}{\section}
\newcommand{\SubSection}{\subsection}

\newtheorem{theorem}{Theorem}

\newcommand{\bbeta}{\mbox{\boldmath  $\beta$}}

\newcommand{\bPsi}{\mbox{\boldmath   $\Psi$}}

\newcommand{\blambda}{\mbox{\boldmath $\lambda$}}

\newcommand{\bzeta}{\mbox{\boldmath $\zeta$}}

\def\urltilda{\kern -.15em\lower .7ex\hbox{\~{}}\kern .04em}
\def\urldot{\kern -.10em.\kern -.10em}
\def\urlhttp{http\kern -.10em\lower -.1ex\hbox{:}\kern -.12em\lower 0ex\hbox{/}\kern -.18em\lower 0ex\hbox{/}}

\begin{document}

\title{An Augmented Lagrangian Approach to the Constrained Optimization Formulation of Imaging Inverse Problems }

\author{Manya V. Afonso, Jos\'{e} M. Bioucas-Dias, and M\'{a}rio A. T. Figueiredo \thanks{The authors are with the {\it Instituto de Telecomunica\c{c}\~{o}es} and the
Department of Electrical and Computer Engineering,
{\it Instituto Superior T\'{e}cnico,} 1049-001 Lisboa, {\bf Portugal}. } \thanks{M. Afonso
is supported by a EU Marie-Curie Fellowship (EST-SIGNAL program: {\scriptsize\tt http://est-signal.i3s.unice.fr}); contract MEST-CT-2005-021175.}
\thanks{An earlier and much shorter version of this work has been presented in \cite{CSALSA_ICASSP}.} }

\maketitle

\begin{abstract}
We propose a new fast algorithm for solving one of the standard
approaches to ill-posed linear inverse problems (IPLIP),
where a (possibly non-smooth) regularizer is minimized
under the constraint that the solution explains the observations
sufficiently well. Although the regularizer and constraint are
usually convex, several particular features of these problems
(huge dimensionality, non-smoothness) preclude the use of
off-the-shelf optimization tools and have stimulated a considerable
amount of research. In this paper, we propose a new efficient
algorithm to handle one class of constrained problems (often
known as basis pursuit denoising) tailored to image recovery
applications. The proposed algorithm, which belongs to the family
of augmented Lagrangian methods, can be used to deal with a variety
of imaging IPLIP, including deconvolution and reconstruction from
compressive observations (such as MRI), using either total-variation
or wavelet-based (or, more generally, frame-based) regularization.
The proposed algorithm is an instance of the so-called
{\it alternating direction method of multipliers}, for which
convergence sufficient conditions are known; we show that these
conditions are satisfied by the proposed algorithm. Experiments on
a set of image restoration and reconstruction benchmark problems
show that the proposed algorithm is a strong contender for the
state-of-the-art.
\end{abstract}

\Section{Introduction}
\SubSection{Problem Formulation}
\label{sec:intro}
Image restoration/reconstruction is one of the earliest and most classical linear
inverse problems in imaging, dating back to the 1960's
\cite{Andrews}. In this class of problems,  a noisy indirect observation ${\bf y}$,
of an original image ${\bf x}$, is modeled as
\[
{\bf y} = {\bf B}{\bf x} + {\bf n},
\]
where ${\bf B}$  is the matrix representation of the direct
operator and ${\bf n}$ is noise. As is common, we adopt
the vector notation for images, where the pixels on an
$M\times N$ image are stacked into ${\bf x}$ an $(NM)$-vector in,
e.g., lexicographic order. In the sequel, we denote by $n=MN$
the number of elements of ${\bf x}$, thus ${\bf x}\in\mathbb{R}^n$,
while  ${\bf y}\in\mathbb{R}^m$ ($m$ and $n$ may be different).

In the particular case of image deconvolution, ${\bf B}$
is the matrix representation of a convolution operator. This type of model describes well several
physical mechanisms, such as relative motion between the camera and
the subject (motion blur), bad focusing (defocusing blur), or a
number of other mechanisms \cite{Bertero}.

In more general image reconstruction problems, ${\bf B}$ represents
some linear direct operator, such as tomographic projections
(Radon transform), a partially observed (e.g., Fourier) transform, or the loss
of part of the image pixels.

The problem of estimating ${\bf x}$ from ${\bf y}$ is called a {\it linear inverse
problem} (LIP); for most scenarios of practical interest, this is an {\it ill-posed LIP}
(IPLIP), {\it i.e.}, matrix ${\bf B}$ is singular and/or very ill-conditioned.
Consequently, this IPLIP requires some
sort of regularization (or prior information, in Bayesian
inference terms). One way to regularize the problem of estimating ${\bf x}$, given ${\bf y}$,
consists in a constrained optimization problem of the form
\begin{equation}\label{genconstrained}
\min_{\bf x} \phi(\bf x)\hspace{0.5cm} \mbox{subject to} \hspace{0.5cm} \|{\bf B}{\bf x}-{\bf y}\|_{2} \leq \varepsilon,
\end{equation}where $\phi: \mathbb{R}^n \rightarrow \bar{\mathbb{R}} = \mathbb{R}\cup\{-\infty,+\infty\}$
is the {\it regularizer} or {\it regularization function},
and $\varepsilon \geq 0$ a parameter which depends on the noise variance.
In the case where $\phi({\bf x}) = \|{\bf x}\|_1$, the above problem is usually known
as {\it basis pursuit denoising} (BPD) \cite{ChenDonohoSaunders}. The
so-called {\it basis pursuit} (BP) problem is the particular case of (\ref{genconstrained})
for $\varepsilon = 0$.

In recent years, an explosion of interest in problems
of the form  (\ref{genconstrained}) was sparked by the emergence of
{\it compressive sensing} (CS) \cite{Candes}, \cite{donoho1}. The theory of
CS  provides conditions (on matrix ${\bf B}$ and the
degree of sparseness of the original ${\bf x}$) under which a solution
of (\ref{genconstrained}), for $\phi({\bf x}) = \|{\bf x}\|_1$,
 is an optimal (in some sense) approximation
to the ``true" ${\bf x}$.

In most signal/image recovery and CS problems, the best results are obtained
with non-smooth regularizers; typical examples are the {\it total variation}
(TV) \cite{Candes}, \cite{ROF} and  $\ell_1$ ($\phi({\bf x}) = \|{\bf x}\|_1$) norms.

\subsection{Analysis and Synthesis Formulations}
\label{sec:synth_analysis}
In a frame-based representation, the unknown image ${\bf x}$ can be represented as a linear
combination of the elements of some frame, {\it i.e.},
${\bf x} = {\bf W}\bbeta$, where $\bbeta\in \mathbb{R}^d$, and the
columns of the $n\times d$ matrix ${\bf W}$ are the elements of a
wavelet\footnote{We will use the generic term ``wavelet" to mean any wavelet-like
multi-scale representation, such as ``curvelets", ``beamlets", ``ridgelets".}
frame (an orthogonal basis or an overcomplete dictionary).
The coefficients of this representation are then estimated
from the noisy image, under one of the well-known
sparsity inducing regularizers, such as the $\ell_1$ norm
(see \cite{GEM}, \cite{DaubechiesDefriseDeMol}, \cite{EladCVPR2006},
\cite{FigueiredoDiasNowak2007}, \cite{FigueiredoNowak2003}, and further references therein).
Formally, this leads to the following constrained optimization problem:
\begin{equation}
\widehat{\bbeta} = \arg\min_{\bbeta} \phi (\bbeta)\hspace{0.5cm} \mbox{subject to}\hspace{0.5cm} \|{\bf B W}\bbeta-{\bf y}\|_{2} \leq \varepsilon.
\label{problem_synthesis}
\end{equation}
This formulation is referred to as the {\it synthesis approach}
\cite{EladMilanfarRubinstein}, \cite{SelesnickFigueiredo}, since
it is based on a synthesis equation: ${\bf x}$ is
{\it synthesized} from its representation coefficients
(${\bf x} = {\bf W}\bbeta$) which are the object of the
estimation criterion. Naturally, the estimate of ${\bf x}$ is
$\widehat{\bf x} = {\bf W}\widehat{\bbeta}$. Of course,
(\ref{problem_synthesis}) has the form (\ref{genconstrained}) with
${\bf B W}$ replacing ${\bf B}$.

An alternative formulation applies a regularizer directly
to the unknown image, leading to criteria of the form
(\ref{genconstrained}), usually called {\it analysis
approaches}, since they are based on a regularizer that
analyzes the image itself, rather than
the coefficients of a representation thereof.
Arguably, the best known and most often used
regularizer in analysis approaches to image restoration is
the total variation (TV) norm \cite{Chan_et_al_2005}, \cite{ROF}.

Wavelet-based analysis approaches are also possible and
have the form
\begin{equation}
\label{eq:bp3_analysis}
\min_{{\bf x}} \phi({\bf P} {\bf x})\hspace{0.5cm} \mbox{subject to}
\hspace{0.5cm} \| {\bf B\, x}- {\bf y}\|_2^2 \leq
\varepsilon,
\end{equation}
where ${\bf P}$ is some linear operator (a matrix)
\cite{EladMilanfarRubinstein}. In this paper, we always assume
that ${\bf P}$ is the analysis operator associated with 1-tight
(Parseval) frame, thus ${\bf P}^H{\bf P} = {\bf I}$ \cite{Mallat}.

\SubSection{Previous Algorithms}
\label{sec:previous}
If the regularizers are convex, problems (\ref{genconstrained})--(\ref{eq:bp3_analysis})
are convex, but the very high dimension (at least $\geq 10^4$, often $\gg 10^5$)
of ${\bf x}$ and ${\bf y}$ precludes the direct application of
off-the-shelf optimization algorithms.
This difficulty is further amplified by the following fact:
for any problem of non-trivial dimension, matrices
${\bf B}$, ${\bf W}$, or ${\bf P}$
cannot be stored explicitly, and it is costly, even impractical,
to access portions (lines, columns, blocks) of them.
On the other hand, matrix-vector products involving these
matrices (or their conjugate transposes, denote by $(\cdot)^{H}$)
can be computed quite efficiently. For example, if the columns of
${\bf W}$ contain a wavelet basis or a tight wavelet frame,
any multiplication of the form ${\bf W}{\bf v}$ or ${\bf W}^H{\bf v}$
can be performed by a fast  transform algorithm \cite{Mallat}.
Similarly, if ${\bf B}$ represents a convolution, products by
${\bf B}$ or ${\bf B}^H$ can be performed with the help
of the fast Fourier transform (FFT). These facts have
stimulated the development of special purpose
methods, in which the only operations involving matrices
are matrix-vector products.

Most state-of-the-art methods for dealing with linear inverse
problems, under convex, non-smooth regularizers (namely, TV and $\ell_1$),
consider, rather than (\ref{genconstrained}), the unconstrained problem
\begin{equation}
\min_{{\bf x}} \frac{1}{2}\| {\bf B\, x} -
{\bf y}\|_2^2 + \tau \, \phi ({\bf x}),\label{unconstrained}
\end{equation}
where $\tau \in \mathbb{R}_+$ is the so-called regularization parameter.
Of course, problems (\ref{genconstrained}) and (\ref{unconstrained})
are equivalent, in the following sense: for any $\varepsilon > 0$ such that
problem (\ref{genconstrained}) is feasible, a solution of (\ref{genconstrained})
is either the null vector, or else is a solution of
(\ref{unconstrained}), for some $\tau > 0$  \cite{GPSR}, \cite{Roc70}.
For solving problems of the form (\ref{unconstrained}), some of the
state-of-the-art algorithms belong to the iterative shrinkage/thresholding (IST)
family \cite{CombettesSIAM, DaubechiesDefriseDeMol, FigueiredoNowak2003, Hale},
and its two-step  (TwIST \cite{TwIST} and FISTA \cite{FISTA}) and accelerated (SpaRSA \cite{SpaRSA_SP}) variants.
These methods were shown to be considerably faster than earlier methods, including $l1\verb=_=ls$
\cite{KimKLBG07} and the codes in the $\ell_1$-{\it magic}\footnote{Available at {\tt http://www.l1-magic.org}}
and  the {\it SparseLab}\footnote{Available at {\tt http://sparselab.stanford.edu}}
toolboxes.

A key ingredient of most of these algorithms is the so-called shrinkage/thresholding/denoising
function, which is the Moreau proximal mapping  of the regularizer $\phi$ \cite{CombettesSIAM}.
Formally, this function $\bPsi_{\tau\phi}:\mathbb{R}^n\rightarrow \mathbb{R}^n$ is defined as
\begin{equation}
\bPsi_{\tau\phi}({\bf y}) = \arg\min_{{\bf x}} \frac{1}{2}\|{\bf x}-{\bf y}\|_2^2 + \tau\phi({\bf x}).
\label{MPM}
\end{equation}
Notice that if $\phi$ is proper and convex, the function being minimized is
proper and strictly convex, thus the minimizer exists and is unique making the
function well defined \cite{CombettesSIAM}. For some choices of $\phi$, the
corresponding $\bPsi_{\tau\phi}$ have well known closed forms. For example,
if $\phi({\bf x}) \equiv \|{\bf x}\|_1$, then $\bPsi_{\tau\phi}({\bf y}) = \mbox{soft}({\bf y}, \tau)$, where
$\mbox{soft}(\cdot, \tau)$ denotes the component-wise application of the
{\it soft-threshold} function $y \mapsto \mbox{sign}(y)\max\{|y|-\tau,0\}$.

In \cite{SALSA_SSP},\cite{SALSA_TIP}, we proposed a new algorithm called {\it split augmented Lagrangian
shrinkage algorithm} (SALSA), to solve unconstrained optimization problems of the
form (\ref{unconstrained}) based on variable splitting \cite{Courant}, \cite{Wang}.
The idea is to transform the unconstrained problem (\ref{unconstrained}) into
a constrained one via a variable splitting ``trick", and then attack this
constrained problem using an {\it augmented Lagrangian} (AL) method \cite{NocedalWright}.
AL is known to be equivalent to the Bregman iterations recently proposed to
handle imaging inverse problems (see \cite{YinOsherGoldfarbDarbon} and references
therein). We prefer the AL perspective, rather than the Bregman iterative view,
as it is a more standard and elementary tool (covered in most optimization
textbooks). On several benchmark experiments (namely image deconvolution,
recovery of missing pixels, and reconstruction from partial Fourier observations)
using either frame-based or TV-based regularization, SALSA was found to be
faster than the previous state-of-the-art methods
FISTA \cite{FISTA}, TwIST \cite{TwIST}, and SpaRSA \cite{SpaRSA_SP}.

Although it is usually easier to solve an
unconstrained problem than a constrained one,
formulation (\ref{genconstrained}) has an important
advantage: parameter $\varepsilon$ has a clear
meaning (it is proportional to the noise standard deviation) and is much
easier to set than parameter $\tau$ in (\ref{unconstrained}).
Of course, one may solve (\ref{genconstrained}) by solving
(\ref{unconstrained}) and searching for the ``correct" value of
$\tau$ that makes (\ref{unconstrained}) equivalent to
(\ref{genconstrained}). Clearly, this is not efficient,
as it involves solving many instances of (\ref{unconstrained}).
Obtaining fast algorithms for solving (\ref{genconstrained}) is thus
an important research front.

There are few efficient algorithms to solve (\ref{genconstrained})-(\ref{eq:bp3_analysis}) in an image recovery context: ${\bf x}$ and ${\bf y}$ of dimension $\geq 10^4$ (often $\geq 10^6$), ${\bf B}$ representing an operator, and $\phi$ a convex, non-smooth function. A notable exception is the recent SPGL1 \cite{SPGL1}, which (as its name implies) is specifically designed for $\ell_1$ regularization. As shown in \cite{SPGL1}, the methods for solving (\ref{genconstrained}) available in the $\ell_1$-magic package are quite inefficient for large problems. General purpose  methods, such as the SeDuMi package\footnote{Available at {\tt http://sedumi.ie.lehigh.edu}}, are simply not applicable when ${\bf B}$ is not an actual matrix, but an operator.

Another efficient algorithm for solving problems of the form
(\ref{genconstrained}) is the recently proposed NESTA \cite{NESTA},
which is based on Nesterov's first-order methods which achieve an
optimal convergence rate by coupling smoothing techniques with an
improved gradient method \cite{Nesterov}, \cite{Nesterov2005}. NESTA allows
for minimizing either the $\ell_1$ or TV norm, and also allows using
synthesis or analysis formulations. Nesterov's first-order method was
also recently adopted in \cite{Dahl} to perform TV-regularized image
denoising, deblurring, and inpainting.

The Bregman iterative algorithm (BIA) was recently proposed to solve
(\ref{genconstrained}) with $\varepsilon = 0$,  but is not
directly applicable when $\varepsilon > 0$ \cite{YinOsherGoldfarbDarbon}.
To deal with the case of $\varepsilon > 0$, it was suggested
that the BIA for $\varepsilon = 0$ is used and stopped when
$\|{\bf B}{\bf x}-{\bf y}\|_{2}\leq \varepsilon$ \cite{CaiOsherZhen},
\cite{YinOsherGoldfarbDarbon}. Clearly, that approach is not
guaranteed to find a good solution, since it depends strongly on the initialization; {\it e.g.},
if the algorithm starts at a feasible point, it will immediately stop,
although the point may be far from a minimizer of $\phi$.

\SubSection{Proposed Approach}


In this paper, we introduce an algorithm for solving
optimization problems of the form (\ref{genconstrained}). The original
constrained problem (\ref{genconstrained}) is transformed
into an unconstrained one by adding the indicator function of
the feasible set, the ellipsoid
$\{ {\bf x}:\|{\bf B x}-{\bf y}\| \leq \varepsilon\}$, to the objective
in (\ref{genconstrained}). The resulting unconstrained problem
is then transformed into a different constrained problem, by the
application of a variable splitting operation; finally, the
obtained constrained problem is dealt with using
the {\it alternating direction method of multipliers} (ADMM)
\cite{EcksteinBertsekas}, \cite{Gabay}, \cite{Glowinski},
which belongs to the family of {\it augmented Lagrangian} (AL)
techniques \cite{NocedalWright}. Since (as SALSA), the proposed
method uses variable splitting and AL optimization, we call it
C-SALSA (for {\it constrained}-SALSA).


The resulting algorithm is more general than SPGL1, in
the sense that it can be used with any convex regularizer
$\phi$ for which the corresponding Moreau proximity
 operator (see \cite{CombettesSIAM}), has closed form
 or can be efficiently computed. In this paper, we
 will show examples of C-SALSA where ${\bf x}$ is an
  image, $\phi$ is the TV norm \cite{ROF},
and $\bPsi_{\tau\phi}$ is computed using Chambolle's
algorithm \cite{Chambolle}. Another classical choice
 which we will demonstrate is the $\ell_1$ norm, which leads to $\bPsi_{\tau\phi}({\bf y}) = \mbox{soft}({\bf y}, \tau)$.

C-SALSA is experimentally shown to efficiently solve
image recovery problems, such as MRI reconstruction from CS-type partial Fourier
observations using TV regularization, and image deblurring
using wavelet-based or TV regularization, faster than SPGL1
and NESTA.

\subsection{Organization of the Paper}
The paper is organized as follows. Section \ref{sec:tools}
describes the basic ingredients of C-SALSA: variable splitting,
augmented Lagrangians, and the ADMM. Section \ref{sec:salsa}
contains the derivation leading to C-SALSA.
Section \ref{sec:experiments} reports experimental results, and
Section \ref{sec:conclusions} ends the paper with a few remarks and
pointers to future work.

\Section{Basic Ingredients}
\label{sec:tools}

\subsection{Variable Splitting}

Consider an unconstrained optimization problem
\begin{equation}
\min_{{\bf u}\in \mathbb{R}^n} f_1({\bf u}) +
f_2\left({\bf G}{\bf u}\right),\label{unconstrained_basic}
\end{equation}
where ${\bf G} \in \mathbb{R}^{d\times n}$, $f_{1}:\mathbb{R}^{n}
\rightarrow \bar{\mathbb{R}}$, and $f_{2}:\mathbb{R}^{d} \rightarrow \bar{\mathbb{R}}$.
Variable splitting (VS) is a  simple procedure that consists in
creating a new variable, say ${\bf v}$,
to serve as the argument of $f_2$, under the constraints that ${\bf v} = {\bf Gu}$, {\it i.e.},
\begin{equation}
 {\displaystyle \min_{{\bf u} \in \mathbb{R}^n ,\, {\bf v}\in\mathbb{R}^d }}\;\;  f_1({\bf u}) + f_2({\bf v}),
 \hspace{0.6cm} \mbox{subject to} \hspace{0.3cm} {\bf v} = {\bf G}{\bf u}.
\label{constrained_basic}
\end{equation}
The rationale behind VS is that it may be easier to solve the constrained problem (\ref{constrained_basic})
than it is to solve its equivalent unconstrained counterpart (\ref{unconstrained_basic}).

VS has been recently used in several image processing applications.
A VS method was used in \cite{Wang} to obtain a fast
algorithm for TV-based restoration. Variable splitting was also
used in \cite{BioucasFigueiredo2008} to handle problems involving compound
regularizers. In \cite{BioucasFigueiredo2008} and \cite{Wang}, the constrained problem (\ref{constrained_basic}) is attacked using a quadratic penalty approach, i.e., by solving
\begin{equation}
 {\displaystyle \min_{{\bf u}\in \mathbb{R}^n,\, {\bf v}\in\mathbb{R}^d}}\;\;
   f_1({\bf u}) + f_2({\bf v}) +
 \frac{\alpha}{2}\, \|{\bf G u} -{\bf v}\|_2^2,
\label{quadratic_penalty}
\end{equation}
by alternating minimization with respect to ${\bf u}$ and ${\bf v}$, while
slowly taking $\alpha$ to very large values (a {\it continuation} process),
to force the solution of (\ref{quadratic_penalty}) to approach that
of (\ref{constrained_basic}), which in turn is equivalent to (\ref{unconstrained_basic}). The rationale of these methods is that each step of this alternating minimization may be much easier than the original
unconstrained problem (\ref{unconstrained_basic}). The drawback is that
as $\alpha$ grows, the intermediate minimization problems
become increasingly ill-conditioned, thus causing numerical problems
\cite{NocedalWright}.

A similar VS approach underlies the recently
proposed split-Bregman methods \cite{GoldsteinOsher}; however,
instead of using a quadratic penalty technique, those methods
attack the constrained problem directly using a Bregman iterative
algorithm \cite{YinOsherGoldfarbDarbon}, which is known to be
equivalent to the augmented Lagrangian method \cite{Esser}, \cite{Setzer},
 \cite{YinOsherGoldfarbDarbon}.

\SubSection{Augmented Lagrangian}
Consider the constrained optimization problem with linear equality constraints
\begin{equation} \min_{{\bf z}\in \mathbb{R}^n}  E({\bf z})\hspace{0.5cm}
 \mbox{subject to}\hspace{0.5cm} {\bf A z - b } = \mbox{\boldmath $0$}, \label{constrained_linear}
\end{equation}
where ${\bf b} \in \mathbb{R}^p$ and ${\bf A}\in \mathbb{R}^{p\times n}$,
{\it i.e.}, there are $p$ linear equality constraints.
The so-called augmented Lagrangian for this problem is
defined as
\begin{equation}
{\cal L}_A ({\bf z},\blambda,\mu) = E({\bf z}) + \blambda^T ({\bf b-Az}) +
\frac{\mu}{2}\,  \|{\bf Az-b}\|_2^2,\label{augmented_L}
\end{equation}
where $\blambda \in \mathbb{R}^p$ is a vector of Lagrange multipliers
and $\mu \geq 0$ is called the AL penalty parameter \cite{NocedalWright}.
The so-called {\it augmented Lagrangian method} (ALM) \cite{NocedalWright},
also known as the {\it method of multipliers} (MM) \cite{Hestenes},
\cite{Powell}, iterates between minimizing ${\cal L}_A ({\bf z},\blambda,\mu)$
with respect to ${\bf z}$, keeping $\blambda$ fixed, and updating
$\blambda$, until some convergence criterion is satisfied:

\vspace{0.3cm}
\begin{algorithm}{ALM/MM}{
\label{alg:salsa1}}
Set $k=0$, choose $\mu > 0$ and  $\blambda_0$.\\
\qrepeat\\
     ${\bf z}_{k+1} \in \arg\min_{{\bf z}} {\cal L}_A ({\bf z},\blambda_k,\mu)$\\
     $\blambda_{k+1} = \blambda_{k} + \mu ({\bf Az}_{k+1} - {\bf b})$\\
     $k \leftarrow k + 1$
\quntil some stopping criterion is satisfied.
\end{algorithm}
\vspace{0.3cm}

It is also possible (and even recommended) to update the value
of $\mu$ in each iteration \cite{Bazaraa}, \cite{NocedalWright}.
However, unlike in the quadratic penalty approach, the ALM/MM does
not require $\mu$ to be taken to infinity to guarantee convergence
to the solution of the constrained problem (\ref{constrained_linear}).

After a straightforward complete-the-squares procedure,
the terms added to $E({\bf z})$ in the augmented Lagrangian  ${\cal L}_A ({\bf z},\blambda_k,\mu)$ can be written as a single quadratic term (plus a
constant independent of ${\bf z}$, thus irrelevant to the ALM/MM),
leading to the following alternative form of the algorithm (which makes
clear its equivalence with the Bregman iterative method \cite{YinOsherGoldfarbDarbon}):

\vspace{0.3cm}
\begin{algorithm}{ALM/MM (version II)}{
\label{alg:salsa2}}
Set $k=0$, choose $\mu > 0$ and  ${\bf d}_0$.\\
\qrepeat\\
     ${\bf z}_{k+1} \in \arg\min_{{\bf z}} E({\bf z}) + \frac{\mu}{2}\|{\bf Az-d}_k\|_2^2$\\
     ${\bf d}_{k+1} = {\bf d}_{k} - ( {\bf Az}_{k+1} - {\bf b})$\\
     $k \leftarrow k+1$
\quntil some stopping criterion is satisfied.
\end{algorithm}
\vspace{0.3cm}

\SubSection{ALM/MM for Variable Splitting and ADMM}
The constrained problem (\ref{constrained_basic}) can be written as (\ref{constrained_linear}) by
defining $E({\bf z})  \equiv f_1({\bf u}) + f_2({\bf v})$ and setting
\begin{equation}
{\bf z} \equiv  \left[{\bf u}^T \; {\bf v}^T\right]^T,\hspace{0.5cm}
{\bf b} = {\bf 0}, \hspace{0.5cm}
{\bf A} =  [{\bf G} \;\; -{\bf I}].
\end{equation}
With these definitions in place, Steps 3 and 4 of the ALM/MM (version II)
become
\begin{align}
\left({\bf u}_{k+1} , {\bf v}_{k+1} \right)  & \in  \arg\min_{{\bf u},{\bf v}} f_{1}({\bf u})
 + f_{2}({\bf v})  + \frac{\mu}{2} \|{\bf G}{\bf u} - {\bf v} - {\bf d_k}\|_2^2 \nonumber\\
{\bf d}_{k+1}  & = {\bf d}_{k} - ({\bf G}{\bf u}_{k+1} - {\bf v}_{k+1}).
\end{align}
The minimization problem yielding $\left({\bf u}_{k+1} , {\bf
v}_{k+1} \right)$ is not trivial since, in general, it involves a
non-separable quadratic term and possibly non-smooth terms. A
natural approach is to use a non-linear block-Gauss-Seidel (NLBGS)
technique which alternates between minimizing with respect to ${\bf
u}$ and ${\bf v}$ while keeping the other fixed. Of course this
raises several questions: for a given ${\bf d}_k$, how much
computational effort should be spent in this problem? Does the NLBGS
procedure converge? The simplest answer to these questions is given
in the form of the so-called {\it alternating direction method of
multipliers} (ADMM) \cite{EcksteinBertsekas}, \cite{Gabay},
\cite{Glowinski}, which is simply an ALM/MM in which only one NLBGS
step is performed in each outer iteration.

\vspace{0.3cm}
\begin{algorithm}{ADMM}{
\label{alg:salsa3}}
Set $k=0$, choose $\mu > 0$, ${\bf v}_0$, ${\bf d}_{0}$.\\
\qrepeat\\
$  {\bf u}_{k+1}  \in  \arg\min_{{\bf u}} f_{1}({\bf u})
 + \frac{\mu}{2} \|{\bf G}{\bf u} - {\bf v}_{k} - {\bf d}_{k}\|_2^2$\\
$  {\bf v}_{k+1}  \in  \arg\min_{{\bf v}} f_{2}({\bf v})
 + \frac{\mu}{2} \|{\bf G}{\bf u}_{k+1} - {\bf v} - {\bf d}_{k}\|_2^2$\\
     ${\bf d}_{k+1} = {\bf d}_{k} - ({\bf Gu}_{k+1} - {\bf v}_{k+1})$\\
     $k \leftarrow k+1$
\quntil some stopping criterion is satisfied.
\end{algorithm}
\vspace{0.3cm}

For later reference, we now recall the theorem by Eckstein and Bertsekas
\cite[Theorem 8]{EcksteinBertsekas} in which convergence of (a generalized version of) ADMM is
shown.
\vspace{0.3cm}
\begin{theorem}[Eckstein-Bertsekas, \cite{EcksteinBertsekas}]
\label{th:Eckstein}{\sl Consider problem  (\ref{unconstrained_basic}), where
 ${\bf G}$ has full column rank, and
 $f_1$ and $f_2$ are closed, proper, convex functions. Consider arbitrary
 $\mu>0$ and ${\bf u}_0, {\bf d}_0, {\bf v}_0\in \mathbb{R}^p$.
 Let $\{\eta_k \geq 0, \; k=0,1,...\}$ and $\{\nu_k \geq 0,
\; k=0,1,...\}$ be two sequences such that
\[
\sum_{k=0}^\infty \eta_k < \infty \;\;\;\mbox{and} \;\;\; \sum_{k=0}^\infty \nu_k <
\infty.
\]
Consider three sequences $\{{\bf u}_k \in \mathbb{R}^{n}, \; k=0,1,...\}$, $\{{\bf v}_k
\in \mathbb{R}^{d}, \; k=0,1,...\}$, and $\{{\bf d}_k \in \mathbb{R}^{d}, \; k=0,1,...\}$
that satisfy
\begin{eqnarray}
\left\| {\bf u}_{k+1} - \arg\min_{{\bf u}} f_{1}({\bf u})
 + \frac{\mu}{2} \|{\bf G}{\bf u} \! - \!{\bf v}_k \! -\! {\bf d}_k\|_2^2 \right\| & \leq & \eta_k   \nonumber\\
\left\| {\bf v}_{k+1}  - \arg\min_{{\bf v}} f_{2}({\bf v})
 + \frac{\mu}{2} \|{\bf G}{\bf u}_{k+1} \! - \! {\bf v} \! - \! {\bf d}_k\|_2^2 \right\| & \leq & \nu_k  \nonumber\\
 {\bf d}_{k+1}  = {\bf d}_{k} - ({\bf G\, u}_{k+1} - {\bf v}_{k+1}).\nonumber
 \end{eqnarray}
Then, if (\ref{unconstrained_basic}) has a solution, say ${\bf u}^*$, then
the sequence $\{{\bf u}_k\}$  converges to ${\bf u}^*$.
If (\ref{unconstrained_basic}) does not have a solution, then at least
one of the sequences $\{{\bf u}_k \}$ or $\{{\bf d}_k\}$ diverges.}
\end{theorem}
\vspace{0.3cm}

Notice that the ADMM as defined above (if each step is
implemented exactly) generates sequences
$\{{\bf u}_k\}$, $\{{\bf v}_k\}$, and $\{{\bf d}_k \}$ that
satisfy the conditions in Theorem \ref{th:Eckstein} in a
strict sense ({\it i.e.}, with $\eta_k = \nu_k = 0$).
The remaining key condition for convergence is then
that ${\bf G}$ has full column rank. One of the important
corollaries of this theorem is that it is not necessary
to exactly solve the minimizations in lines 3 and 4 of ADMM;
as long as the sequence of errors are absolutely summable,
convergence is not compromised.

The proof of Theorem \ref{th:Eckstein} is based on the
equivalence between ADMM and the Douglas-Rachford Splitting (DRS)
applied to the dual of  problem (\ref{unconstrained_basic}).
The DRS was recently used for image recovery problems in
\cite{CombettesPesquet}. For recent and comprehensive reviews
of ALM, ADMM, DRS, and their relationship with Bregman and
split-Bregman methods, see \cite{Esser}, \cite{Setzer}.

\subsection{A Variant of ADMM}\label{sec:variant}

Consider a generalization
of problem (\ref{unconstrained_basic}), where instead of two functions,
there are $J$ functions, that is,
\begin{equation}
\min_{{\bf u}\in \mathbb{R}^d} \; \sum_{j=1}^J g_j ({\bf H}^{(j)}\,{\bf u}),\label{unconstrained_compound}
\end{equation}
where $g_j : \mathbb{R}^{p_j} \rightarrow \bar{\mathbb{R}}$ are closed, proper,
convex functions, and ${\bf H}^{(j)} \in \mathbb{R}^{p_j \times d}$ are arbitrary
matrices. The minimization problem (\ref{unconstrained_compound}) can be written
as (\ref{unconstrained_basic}) using the following correspondences: $f_1 = 0$,
\begin{equation}
{\bf G}  =  \left[ \begin{array}{c} {\bf H}^{(1)} \\ \vdots \\
{\bf H}^{(J)} \end{array}\right] \in \mathbb{R}^{p\times d},\label{eq:stackedG}
\end{equation}
where $p=p_1 + \dots + p_J$, and $f_2 : \mathbb{R}^{p\times d} \rightarrow \bar{\mathbb{R}}$ given by
\begin{equation}
f_2({\bf v}) = \sum_{j=1}^J g_j ({\bf v}^{(j)}),\label{eq:f_2}
\end{equation}
where ${\bf v}^{(j)} \in \mathbb{R}^{p_j}$ and ${\bf v} = [({\bf v}^{(1)})^T,\dots,({\bf v}^{(J)})^T]^T \in \mathbb{R}^{p}$.
We now simply apply ADMM (as given in the previous subsection), with
 \[
{\bf d}_k = \left[ \begin{array}{c}{\bf d}^{(1)}_k\\ \vdots \\{\bf d}^{(J)}_k\end{array}\right], \;\;\;
{\bf v}_k = \left[ \begin{array}{c}{\bf v}^{(1)}_k\\ \vdots \\{\bf v}^{(J)}_k\end{array}\right].
\]
Moreover, the fact that $f_1 = 0$ turns Step 3 of the algorithm into a
simple quadratic minimization problem, which has a unique solution if ${\bf G}$
has full column rank:
\begin{eqnarray}
\arg\min_{{\bf u}} \; \bigl\|{\bf G\, u} - \bzeta_k\bigr\|_2^2 & = & \left({\bf G}^H {\bf G}\right)^{-1} {\bf G}^H \bzeta_k,
\label{eq:quadratic_problem}\\
& & \hspace{-2cm} = \; \; \biggl[ \sum_{j=1}^J
({\bf H}^{(j)})^H{\bf H}^{(j)}\biggr]^{-1} \sum_{j=1}^J \bigl({\bf H}^{(j)}\bigr)^H \bzeta^{(j)}_k,\nonumber
\end{eqnarray}
where $\bzeta_k = {\bf v}_k + {\bf d}_k$ (and, naturally, $\bzeta_k^{(j)} = {\bf u}_k^{(j)} + {\bf d}_k^{(j)}$)
and the second equality results from the particular structure of ${\bf G}$ in (\ref{eq:stackedG}).

Furthermore, our particular way of mapping problem (\ref{unconstrained_compound}) into
problem (\ref{unconstrained_basic}) allows decoupling the minimization in Step 4
of ADMM into a set of $J$ independent ones. In fact,
\begin{equation}
 {\bf v}_{k+1} \leftarrow  \arg\min_{{\bf v}} f_{2}({\bf v})
 + \frac{\mu}{2} \|{\bf G\, u}_{k+1} - {\bf v} - {\bf d}_k\|_2^2 \label{update_u}
 \end{equation}
 can be written as
\begin{gather}
\left[ \begin{array}{l} {\bf v}^{(1)}_{k+1}\\ \vdots \\ {\bf v}^{(J)}_{k+1}
  \end{array}\right]
 \leftarrow   \arg\min_{{\bf v}^{(1)}, \dots , {\bf v}^{(J)}} \; g_1 ({\bf v}^{(1)}) + \dots + g_J ({\bf v}^{(J)}) + \nonumber \\
 +
\frac{\mu}{2}\left\| \left[ \begin{array}{c} {\bf H}^{(1)} \\ \vdots \\ {\bf H}^{(J)} \end{array}\right]  {\bf u}_{k+1} -
\left[ \begin{array}{c} {\bf v}^{(1)} \\ \vdots \\ {\bf v}^{(J)}\end{array}\right] - \left[\begin{array}{l} {\bf d}^{(1)}_k \\
\vdots \\ {\bf d}^{(J)}_k\end{array}\right]  \right\|_2^2.\nonumber
\end{gather}
Clearly, the minimizations with respect to ${\bf u}^{(1)},\dots,{\bf u}^{(J)}$
are decoupled, thus can be solved separately, leading to
\begin{equation}
 {\bf v}^{(j)}_{k+1}  \leftarrow  \arg\min_{{\bf v}\in \mathbb{R}^{p_j}} \;  g_j ({\bf v}) +
\frac{\mu}{2} \, \bigl\| {\bf v}  - {\bf s}^{(1)}_k\bigr\|_2^2 , \label{eq:Moreau_1}
\end{equation}
for $j=1,...,J$, where
\[
{\bf s}^{(j)}_k =  {\bf H}^{(j)} {\bf u}_{k+1} - {\bf d}^{(j)}_k.
\]

Since this algorithm is exactly an ADMM, and since all the functions
$g_j$, for $j=1,...,J$, are closed, proper, and convex, convergence is
guaranteed if ${\bf G}$ has full column rank. Actually, this full column
rank condition is also required for the inverse in (\ref{eq:quadratic_problem}) to exist.
Finally, notice that the update equations in  (\ref{eq:Moreau_1})
can be written as
\begin{equation}
{\bf v}_{k+1}^{(j)} = \bPsi_{g_j/\mu}( {\bf s}^{(j)}_k ).
\end{equation}
where the $\bPsi_{g_j/\mu}$ are, by definition, the {\it Moreau proximal mappings}
of $g_1/\mu,...,g_J/\mu$.

In summary, the variant of ADMM (herein referred to as ADMM-2) that results
from the formulation just presented is described in the following algorithmic
framework.

\vspace{0.3cm}
\begin{algorithm}{ADMM-2}{
\label{alg:admm-2}}
Set $k=0$, choose $\mu > 0$, ${\bf v}_0^{(1)}$, ..., ${\bf v}_0^{(J)}$, ${\bf d}_{0}^{(1)}$, ..., ${\bf d}_{0}^{(J)}$.\\
\qrepeat\\
\qfor $i=1,...,J$\\
\qdo $ \bzeta_k^{(i)} = {\bf v}_k^{(i)} + {\bf d}_k^{(i)}$
\qrof\\

$  {\bf u}_{k+1}  = \biggl[ \sum_{j=1}^J
({\bf H}^{(j)})^H{\bf H}^{(j)}\biggr]^{-1} \sum_{j=1}^J \bigl({\bf H}^{(j)}\bigr)^H \bzeta^{(j)}_k $\\
\qfor $i=1,...,J$\\
\qdo $  {\bf v}_{k+1}^{(i)} =  \bPsi_{g_i/\mu} \left(   {\bf H}^{(i)} {\bf u}_{k+1} - {\bf d}^{(i)}_k \right)$\\
  ${\bf d}_{k+1}^{(i)} = {\bf d}_{k}^{(i)} - {\bf H}^{(i)}{\bf u}_{k+1} + {\bf v}_{k+1}^{(i)}$
\qrof\\
     $k \leftarrow k+1$
\quntil some stopping criterion is satisfied.
\end{algorithm}
\vspace{0.3cm}

\Section{Proposed Method}
\label{sec:salsa}

We now apply the algorithmic framework
described in the previous section to the basic problem
(\ref{genconstrained}) (which includes (\ref{problem_synthesis}) as
a special case), as well as the analysis formulation (\ref{eq:bp3_analysis}).

\subsection{Problem (\ref{genconstrained})}
For the constrained optimization problem (\ref{genconstrained}),
the feasible set is the ellipsoid
\begin{equation}
E(\varepsilon,{\bf B},{\bf y}) =
\{{\bf x}\in \mathbb{R}^n:\| {\bf B\, x}- {\bf y}\|_2 \leq \varepsilon\},
\end{equation}
which is possibly infinite in some directions (since ${\bf B}$ may be singular).
Problem (\ref{genconstrained}) can be written as an unconstrained  problem,
with a discontinuous objective,
\begin{equation}
\min_{{\bf x}} \; \phi({\bf x}) + \iota_{E(\varepsilon,{\bf I},{\bf y})} ({\bf B x }) ,
\label{unconstrained_reformulation}
\end{equation}
where $\iota_S:\mathbb{R}^m\rightarrow \bar{\mathbb{R}}$
denotes the indicator function of set $S\subset \mathbb{R}^m$,
\begin{equation}
\iota_S({\bf s})=\left\{
\begin{array}{ll}
0, & \text{if } {\bf s} \in S\\
+\infty, & \text{if } {\bf s} \notin S.
\end{array} \right.\label{indicator_ellipsoid}
\end{equation}
Notice that $E(\varepsilon,{\bf I},{\bf y})$ is simply a closed
$\varepsilon$-radius Euclidean ball centered at ${\bf y}$.

Problem (\ref{unconstrained_reformulation}) has the form (\ref{unconstrained_compound})
with $J=2$ and
\begin{eqnarray}
g_1 & \equiv & \phi \label{eq:def_f1}\\
g_2 & \equiv & \iota_{E(\varepsilon,{\bf I},{\bf y})},\label{eq:def_f2}\\
{\bf H}^{(1)} & \equiv & {\bf I}\\
{\bf H}^{(2)} & \equiv & {\bf B}.
\end{eqnarray}

Instantiating ADMM-2 to this particular case requires the definition
of the Moreau proximal maps associated with $g_1 \equiv \phi$ and $g_2 \equiv \iota_{E(\varepsilon,{\bf I},{\bf y})}$. Concerning
$\phi$, the regularizer, we assume that $\bPsi_{\tau\phi}(\cdot)$
(see (\ref{MPM})) can be computed efficiently. This is of course the case
of $\phi({\bf x}) \equiv \|{\bf x}\|_1$, for which $\bPsi_{\tau\phi}$ is simply a soft threshold.
If $\phi$ is the TV norm, we may use one the fast algorithms available to
compute the corresponding denoising function \cite{Chambolle}, \cite{Dahl}.
The Moreau proximal map of $g_2 \equiv \iota_{E(\varepsilon,{\bf I},{\bf y})}$
is defined as
\begin{equation}
\bPsi_{\iota_{E(\varepsilon,{\bf I},{\bf y})}/\mu}({\bf s}) = \arg\min_{{\bf x}}
 \frac{\iota_{E(\varepsilon,{\bf I},{\bf y})} ({\bf x})}{\mu} + \frac{1}{2}\|{\bf x}-{\bf s}\|_2^2,
\end{equation}
which is obviously independent of $\mu$ and is simply the orthogonal
projection of ${\bf s}$ on the closed $\varepsilon$-radius ball centered at ${\bf y}$:
\begin{equation}
\bPsi_{\iota_{E(\varepsilon,{\bf I},{\bf y})}}({\bf s})  =
{\bf y} + \left\{ \begin{array}{ll}
\varepsilon \, \frac{{\bf s}-{\bf y}}{\|{\bf s}-{\bf y}\|_2}, & \text{if }\;  \|{\bf s} - {\bf y}\|_2 > \varepsilon,\\
{\bf s} - {\bf y}, & \text{if }\;  \|{\bf s} - {\bf y}\|_2 \leq \varepsilon.
\end{array} \right.\label{projection_ellipsoid}
\end{equation}

We are now in a position to instantiate ADMM-2 for solving
(\ref{unconstrained_reformulation}) (equivalently (\ref{genconstrained})).
The resulting algorithm, which we call C-SALSA-1, is as follows.

\vspace{0.3cm}
\begin{algorithm}{C-SALSA-1}{
\label{alg:csalsa1}}
Set $k=0$, choose $\mu > 0$, ${\bf v}_0^{(1)}$, ${\bf v}_0^{(2)}$, ${\bf d}_{0}^{(2)}$, ${\bf d}_{0}^{(2)}$.\\
\qrepeat\\
${\bf r}_k = {\bf v}_0^{(1)} + {\bf d}_0^{(1)} + {\bf B}^H \left({\bf v}_0^{(2)} + {\bf d}_0^{(2)}\right) $\\
$  {\bf u}_{k+1}  = \biggl( {\bf I} + {\bf B}^H {\bf B} \biggr)^{-1} {\bf r}_k $\\
 $  {\bf v}_{k+1}^{(1)} =  \bPsi_{\phi/\mu} \left( {\bf u}_{k+1} - {\bf d}^{(1)}_k \right)$\\
 $  {\bf v}_{k+1}^{(2)} =  \bPsi_{\iota_{E(\varepsilon,{\bf I},{\bf y})}} \left(  {\bf B} {\bf u}_{k+1} - {\bf d}^{(2)}_k \right)$\\
${\bf d}_{k+1}^{(1)} = {\bf d}_{k}^{(1)} - {\bf u}_{k+1} + {\bf v}_{k+1}^{(1)}$\\
${\bf d}_{k+1}^{(2)} = {\bf d}_{k}^{(2)} - {\bf B}{\bf u}_{k+1} + {\bf v}_{k+1}^{(2)}$\\
     $k \leftarrow k+1$
\quntil some stopping criterion is satisfied.
\end{algorithm}
\vspace{0.3cm}

The issue of how to efficiently solve the linear system of
equations in line 4 of C-SALSA-1 will be addressed in Subsection~\ref{sec:computingxk}.

Convergence of  C-SALSA-1 is guaranteed by Theorem \ref{th:Eckstein}
since it is an instance of ADMM with
\begin{equation}
{\bf G} \equiv \left[ \begin{array}{c} {\bf I} \\
{\bf B} \end{array}\right],
\end{equation}
which is a full column rank matrix, and both $\phi$
and $\iota_{E(\varepsilon,{\bf I},{\bf y})}$ are
closed, proper, convex functions.

Finally, notice that to apply C-SALSA-1 to problem (\ref{problem_synthesis})
we simply have to replace ${\bf B}$ with ${\bf BW}$.

\subsection{Problem (\ref{eq:bp3_analysis})}
Problem (\ref{eq:bp3_analysis}) can also be written as an
unconstrained  problem
\begin{equation}
\min_{{\bf x}} \; \phi({\bf P} {\bf x}) + \iota_{E(\varepsilon,{\bf I},{\bf y})} ({\bf B x }) ,
\label{unconstrained_reformulation_2}
\end{equation}
which has the form (\ref{unconstrained_compound}) with $J=2$ and
\begin{eqnarray}
g_1 & \equiv & \phi \label{eq:def_f1_b}\\
g_2 & \equiv & \iota_{E(\varepsilon,{\bf I},{\bf y})},\label{eq:def_f2_b}\\
{\bf H}^{(1)} & \equiv & {\bf P}\\
{\bf H}^{(2)} & \equiv & {\bf B}.
\end{eqnarray}

The resulting ADMM algorithm, called C-SALSA-2, is similar to C-SALSA-1, with
only a few minor differences.

\vspace{0.3cm}
\begin{algorithm}{C-SALSA-2}{
\label{alg:csalsa2}}
Set $k=0$, choose $\mu > 0$, ${\bf v}_0^{(1)}$, ${\bf v}_0^{(2)}$, ${\bf d}_{0}^{(2)}$, ${\bf d}_{0}^{(2)}$.\\
\qrepeat\\
${\bf r}_k ={\bf P} \left({\bf v}_{k}^{(1)} + {\bf d}_{k}^{(1)}\right) + {\bf B}^H\left({\bf v}_{k}^{(2)} + {\bf d}_{k}^{(2)}\right)$\\
${\bf u}_{k+1}  \!= \!\biggl( {\bf P}^H{\bf P} + {\bf B}^H {\bf B} \biggr)^{-1} {\bf r}_k  $\\
 $  {\bf v}_{k+1}^{(1)} =  \bPsi_{\phi/\mu} \left( {\bf P}{\bf u}_{k+1} - {\bf d}^{(1)}_k \right)$\\
 $  {\bf v}_{k+1}^{(2)} =  \bPsi_{\iota_{E(\varepsilon,{\bf I},{\bf y})}} \left(  {\bf B} {\bf u}_{k+1} - {\bf d}^{(2)}_k \right)$\\
${\bf d}_{k+1}^{(1)} = {\bf d}_{k}^{(1)} - {\bf P}{\bf u}_{k+1} + {\bf v}_{k+1}^{(1)}$\\
${\bf d}_{k+1}^{(2)} = {\bf d}_{k}^{(2)} - {\bf B}{\bf u}_{k+1} + {\bf v}_{k+1}^{(2)}$\\
     $k \leftarrow k+1$
\quntil some stopping criterion is satisfied.
\end{algorithm}
\vspace{0.3cm}

In this paper, we assume that ${\bf P}$ is the analysis operator
of a 1-tight (Parseval) frame, thus ${\bf P}^H{\bf P}={\bf I}$ and
line 4 of C-SALSA-2 is similar to line 4 of C-SALSA-1:
\begin{equation}
{\bf u}_{k+1} = \left(  {\bf I} + {\bf B}^H {\bf B} \right)^{-1} {\bf r}_k.\label{eq:the_matrix}
\end{equation}
The issue of how to efficiently solve this linear system will be addressed
in the next subsection.

Since both $\phi$ and $\iota_{E(\varepsilon,{\bf I},{\bf y})}$ are
closed, proper, convex functions, convergence of C-SALSA-2
holds (by Theorem \ref{th:Eckstein}) if
\begin{equation}
{\bf G} = \left[ \begin{array}{c} {\bf P} \\
{\bf B} \end{array}\right],
\end{equation}
is a full column rank matrix. This is of course true if
${\bf P}$ is itself a full column rank matrix, which is the
case if ${\bf P}$ is the analysis operator of a tight frame
\cite{Mallat}.

\SubSection{Solving (\ref{eq:the_matrix})}
\label{sec:computingxk}
As mentioned in Subsection~\ref{sec:previous}, in most imaging
problems of interest, it may not be feasible to explicitly
form matrix ${\bf B}$. This might suggest that it is not easy,
or even feasible, to compute the inverse of $\left({\bf I} +
{\bf B}^H {\bf B} \right)$. However, as shown next, in a
number of problems of interest, this inverse can be
computed very efficiently with $O(n\log n)$ cost.

\vspace{0.3cm}
\subsubsection{Deconvolution with Analysis Formulation}
In the case of analysis formulations of the form (\ref{genconstrained})
or (\ref{eq:bp3_analysis}) to image deconvolution problems, matrix
${\bf B}$ represents a 2D convolution. Consequently, matrix ${\bf B}$
can be factorized as ${\bf B = U}^H {\bf D U}$, where ${\bf U}$ is the unitary matrix
(${\bf U}^H = {\bf U}^{-1} $) representing the discrete Fourier transform (DFT)
and ${\bf D}$ is diagonal. Thus,
\begin{equation}
({\bf B}^H{\bf B} + {\bf I})^{-1} =  {\bf U}^H\left( |{\bf D}|^2 + {\bf I}\right)^{-1}{\bf U},\label{eq:Wiener2}
\end{equation}
where $|{\bf D}|^2$ is the  matrix with squared absolute values
of the entries of ${\bf D}$. Since $|{\bf D}|^2 + {\bf I}$ is diagonal,
its inversion cost is $O(n)$. Products by ${\bf U}$  and ${\bf U}^H$ have
$O(n\log n)$ cost, using the FFT algorithm.

\subsubsection{Deconvolution with Synthesis Formulation}
In this case, as seen in Section~\ref{sec:synth_analysis}, we have
${\bf B W}$ instead of ${\bf B}$, and even if ${\bf B}$  is a convolution,
${\bf BW}$ is not diagonalizable by the DFT. To sidestep this difficulty,
we assume that ${\bf W}$ contains a 1-tight (Parseval) frame
({\it i.e.}, ${\bf W\,W}^H = {\bf I}$). Using the Sherman–-Morrison-–Woodbury (SMW) matrix inversion lemma,
\begin{equation}
\left({\bf W}^H\!{\bf B}^H\!{\bf B W} + {\bf I}\right)^{-1}\!  = \;
 {\bf I} - {\bf W}^H\!\underbrace{{\bf B\!}^H \!\left({\bf BB\!}^H +
 {\bf I}\right)^{-1}\!{\bf B}}_{\bf F} {\bf W},
\label{eq:filter_3}
\end{equation}
thus line 4 of C-SALSA-1 and C-SALSA-2 can be written as
\begin{equation}
{\bf u}_{k+1} =  \left({\bf r}_k - {\bf W}^H {\bf F} {\bf W}\; {\bf r}_k \right).
\label{GS_linear_5}
\end{equation}
Since ${\bf B}$ is a convolution, ${\bf B = U}^H {\bf D U}$, thus multiplying by ${\bf F}$ corresponds to
applying an image filter in the Fourier domain
\[
{\bf F} = {\bf U}^H {\bf D^*} \left(|{\bf D}|^2  +  {\bf I}\right)^{-1} {\bf D U},
\]
which has $O(n\log n)$ cost, since all the matrices in
${\bf D^*} \left( |{\bf D}|^2  + {\bf I}\right)^{-1} {\bf D}$ are diagonal
and the products by ${\bf U}$ and ${\bf U}^H$ are carried out via the
FFT. The  cost of (\ref{eq:filter_3}) will thus be either $O(n\log n)$
or the cost of the products by ${\bf W}^H$ and ${\bf W}$.

For most tight frames used in image processing, there are
fast $O(n\log n)$ algorithms to compute the
products by  ${\bf W}^H$ and ${\bf W}$ \cite{Mallat}.
For example, in the case of  translation-invariant wavelet
transforms, these products can be computed using the
undecimated wavelet transform with $O(n\log n)$ total  cost
\cite{lang}. Curvelets also constitute a Parseval frame for
which fast $O(n\log n)$ implementations of the forward and
inverse transform exist \cite{CandesDemanetDonohoYing}.
Yet another example of a redundant Parseval frame is the
complex wavelet transform, which has $O(n)$ computational
cost \cite{Kingsbury}, \cite{Selesnick}. In conclusion,
for a large class of choices of ${\bf W}$, each iteration
of the SALSA algorithm has $O(n\log n)$ cost.

\vspace{0.3cm}
\subsubsection{Missing Pixels: Image Inpainting}
In the analysis prior formulation of this problem, the observation matrix ${\bf B}$ models the loss of some image pixels; it is thus an $m\times n$
binary matrix, with $m<n$, which can be obtained by taking a subset of
rows of an identity matrix. Due to its particular structure, this matrix
satisfies ${\bf B}{\bf B}^H = {\bf I}$. Using this fact together with
the SMW formula leads to
\begin{eqnarray}
\left({\bf B}^H{\bf B} + {\bf I}\right)^{-1} & = &
{\bf I} - {\bf B}^H \left( {\bf I} + {\bf B}{\bf B}^H\right)^{-1}{\bf B}\\
& = & {\bf I} - \frac{1}{2}\, {\bf B}^H{\bf B}.
\label{eq:mask_matrix}
\end{eqnarray}
Since ${\bf B}^H{\bf B}$ is equal to an identity matrix with some
zeros in the diagonal (corresponding to the positions of the
missing observations), the matrix in (\ref{eq:mask_matrix}) is
diagonal with elements either equal to $1$ or $1/2$.
Consequently,  line 4 of C-SALSA-1 and C-SALSA-2 corresponds to
multiplying this diagonal matrix by ${\bf r}_k$, obviously with
$O(n)$ cost.

In the frame-based synthesis formulation, we have ${\bf B W}$ instead of ${\bf B}$. Using  the SMW formula yet again, and the facts that
${\bf B}{\bf B}^H = {\bf I}$ and ${\bf WW}^H = {\bf I}$, we have
\begin{equation}
\left({\bf W}^H{\bf B}^H{\bf B W} + {\bf I}\right)^{-1} =
 {\bf I} -  \frac{1}{2}\, {\bf W}^H {\bf B}^H {\bf B} {\bf W} .
\label{eq:missing_synth}
\end{equation}
As noted in the previous paragraph, ${\bf A}^H{\bf A}$ is equal
to an identity matrix with zeros in the diagonal, {\it i.e.},
a binary mask. Thus, the multiplication by  ${\bf W}^H {\bf A}^H
{\bf A} {\bf W} $ corresponds to synthesizing the image,
multiplying it by this mask, and computing the representation
coefficients of the result. In conclusion, the cost of line 4
of C-SALSA-1 and C-SALSA-2 is again that of the products by
${\bf W}$ and ${\bf W}^H$, usually $O(n\log n)$.

\vspace{0.3cm}
\subsubsection{Partial Fourier Observations (MRI Reconstruction)}

Finally, we consider the case of partial Fourier observations,
which is used to model  MRI  acquisition and has been the focus of
recent interest due to its connection
to compressed sensing \cite{Candes}, \cite{Lustig}.
In the analysis formulation, ${\bf B}={\bf M}{\bf U}$,
where ${\bf M}$ is an $m\times n$ binary
matrix ($m<n$) again, formed by a subset of rows of the identity,
and ${\bf U}$ is the DFT matrix.
Due to its particular structure, matrix ${\bf M}$ satisfies
${\bf M}{\bf M}^H = {\bf I}$; this fact together with the matrix inversion lemma
leads to
\begin{equation}
\left( {\bf B}^H{\bf B} + {\bf I}\right)^{-1} =
{\bf I} - \frac{1}{2}\, {\bf U}^H {\bf M}^H{\bf M}{\bf U},
\label{eq:mask_matrix_MRI}
\end{equation}
where ${\bf M}^H{\bf M}$ is equal to an identity with some
zeros in the diagonal. Consequently,  the cost of line 4
of C-SALSA-1 and C-SALSA-2 is again that of the products by
${\bf U}$ and ${\bf U}^H$, {\it i.e. } $O(n\log n)$ using the FFT.

In the synthesis case, the observation matrix has the form
${\bf M}{\bf U}{\bf W}$. Clearly, the case is again
similar to (\ref{eq:missing_synth}), but with ${\bf UW}$ and
${\bf W}^H{\bf U}^H$ instead of ${\bf W}$ and ${\bf W}^H$,
respectively. Again, the cost of line 4 of C-SALSA-1 and
C-SALSA-2 is $O(n\log n)$, if the FFT is used to compute
the products by ${\bf U}$ and ${\bf U}^H$ and fast frame
transforms are used for the products by ${\bf W}$ and
${\bf W}^H$.


\subsection{Computational Complexity}
As shown in the previous section, the cost of line 4
of C-SALSA-1 and C-SALSA-2 is $O(n\log n)$.
The other lines of the algorithms simply involve:
(a) matrix-vector products involving ${\bf B}$, ${\bf W}$,
${\bf P}$, or their conjugate transposes, which have $O(n\log n)$
cost; (b) vector additions, which have $O(n)$ cost;
and (c) the computation of the Moreau proximal maps
(lines 5 and 6 of C-SALSA-1 and C-SALSA-2).
In the case of the projections on a ball (line 6), it
is clear from (\ref{projection_ellipsoid}) that the
cost is $O(n)$.

Finally, we consider the computational cost of the Moreau
proximal map of the regularizer $\phi$ (line
5 of C-SALSA-1 and C-SALSA-2). In some cases,
this map can be computed exactly in closed form; for
example, if $\phi({\bf x}) \equiv \|{\bf x}\|_1$, then $\bPsi_{\tau \phi}$ is
simply a soft threshold and the cost is $O(n)$. In other
cases, the Moreau proximal map does not have a closed form
solution; for example, if $\phi({\bf x}) \equiv \mbox{TV}({\bf x})$,
the corresponding $\bPsi_{\tau \phi}$ has to be computed using one
of several available iterative algorithms \cite{Chambolle},
\cite{Dahl}. Most of these iterative algorithms can be implemented
with $O(n)$ cost, although with a factor that depends on
the number of iterations. In our implementation of C-SALSA
we use Chambolle's algorithm \cite{Chambolle}.

In summary, for a wide choice of regularizers and frame
representations, the C-SALSA algorithms have $O(n\log n)$
computational complexity.

\Section{Experiments}
\label{sec:experiments}

In this section, we report results of experiments aimed at comparing
the speed of C-SALSA with that of the current state of the art methods
(that are freely available online): SPGL1\footnote{Available at \url{ http://www.cs.ubc.ca/labs/scl/spgl1}} \cite{SPGL1}, and
NESTA\footnote{Available at \url{http://www.acm.caltech.edu/~nesta}}
\cite{NESTA}.

We consider three standard and
often studied imaging inverse problems:
image deconvolution (using both wavelet and TV-based regularization); image restoration from
missing samples (inpainting); image reconstruction from partial Fourier observations,
which (as mentioned above) has been the focus of much recent interest due to its connection
with compressed sensing and the fact that it models MRI acquisition \cite{Lustig}.

All experiments were performed using MATLAB, on a Windows XP desktop computer with an Intel Pentium-IV $3.0$ GHz
processor and $1.5$ GB of  RAM. The number of calls to the
operators ${\bf B}$ and ${\bf B}^H$, the number of iterations,
CPU times, and MSE values presented are the averages values over
10 runs of each experiment. The number of calls reported for each experiment is the average over the 10 instances, with the minimum and maximum indicated in the parentheses. Since the stopping criteria of the
implementations of the available algorithms differ, to compare
the speed of the algorithms in a way that is as independent as
possible from these criteria, the experimental
protocol that we followed was the following: we first run one
of the algorithms with its stopping criterion, and then run
C-SALSA until the constraint in (\ref{genconstrained}) is
satisfied and the MSE of the estimate is below that
obtained by the other algorithms.

The value of $\varepsilon$ in (\ref{genconstrained}) used in all cases was $\sqrt{m+8\sqrt{m}}\sigma$, where $m$ is the number of observations, and $\sigma$ is the noise standard deviation. The parameter $\mu$ was hand-tuned for fastest convergence.

\subsection{Image Deconvolution with wavelets}

We consider five benchmark deblurring problems \cite{FigueiredoNowak2003},
summarized in Table~\ref{decon_problems}, all on the well-known Cameraman image.
The regularizer is $\phi(\bbeta) = \|\bbeta\|_1$, thus $\bPsi_{\tau\phi}$ is an
element-wise soft threshold. We compare C-SALSA against SPGL1 and NESTA in the synthesis case, and against only NESTA in the analysis case, since SPGL1 is hardwired with $\|{\bf x}\|_1$ as the regularizer, and not $\|{\bf P x}\|_1$. Since the restored images are visually indistinguishable from those obtained in \cite{FigueiredoNowak2003}, and the SNR improvements are also very similar, we simply compare the speed of the algorithms, that is, the number of calls to the operators ${\bf B}$ and ${\bf B}^H$, the number of iterations, and the computation time.

\begin{table}[hbt]
\centering
\caption{Details of the image deconvolution experiments.}\label{decon_problems}
\begin{tabular}{|c|l|l|}
  \hline
 Experiment &  blur kernel  \rule[-0.1cm]{0cm}{0.4cm}  & $\sigma^2$ \\ \hline
1  &$9\times 9$ uniform & $0.56^2$ \\
2A & Gaussian & 2\\
2B & Gaussian & 8\\
3A  & $h_{ij} = 1/(1 + i^2 + j^2)$ & 2 \\
3B  & $h_{ij} = 1/(1 + i^2 + j^2)$ & 8 \\
\hline
\end{tabular}
\end{table}

\begin{figure}
\centering
\subfigure[]{
\includegraphics[width=0.3\textwidth]{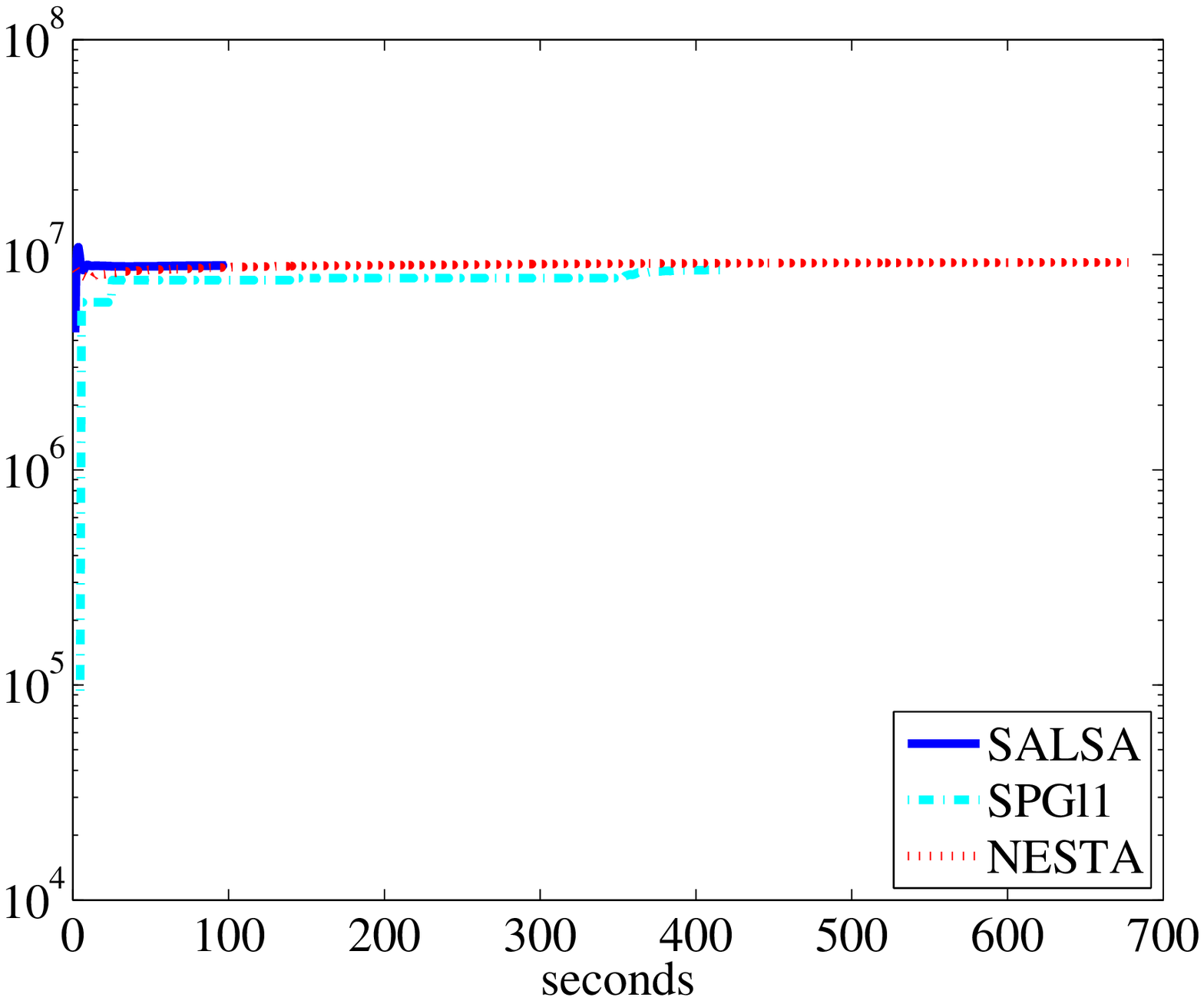}
}
\subfigure[]{
\includegraphics[width=0.3\textwidth]{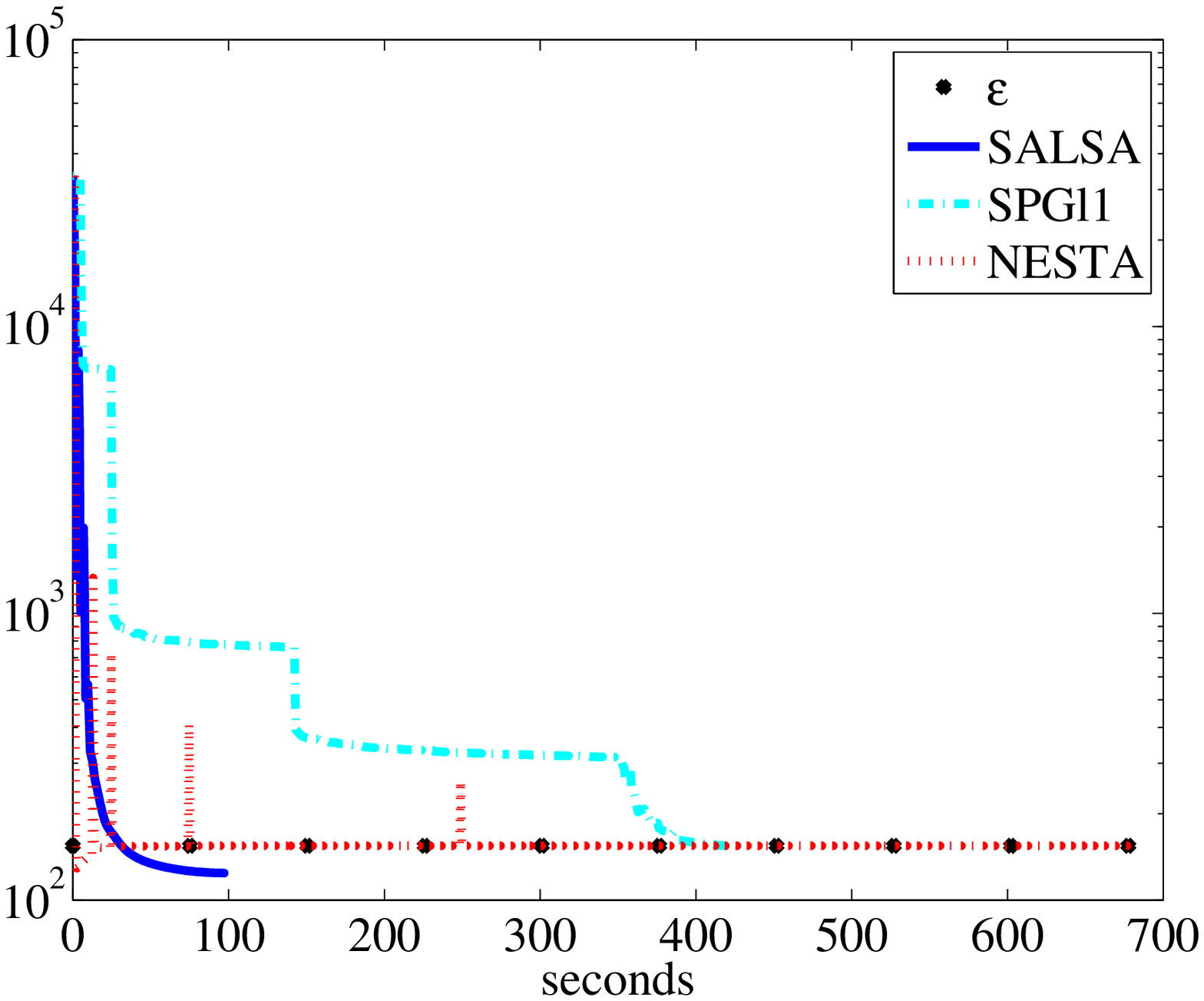}
}
\caption{\footnotesize Image deblurring with wavelets (synthesis prior, redundant Haar wavelets), $9\times 9$ uniform
blur, $\sigma=0.56$: (a) Evolution of the objective function $\|{\bf x}\|_1$ over time; (b) quadratic constraint $\|{\bf AWx}-{\bf y}\|_2$ over time.}
\label{fig:evolutioncriterionl1_redundant}
\end{figure}

In the first set of experiments, ${\bf W}$ is a redundant Haar wavelet frame
with four levels. For the synthesis case, the CPU times taken by each of the algorithms are presented in Table~\ref{tab:deconvl1redundant}. Table~\ref{tab:deconvl1redundant_analysis} presents the corresponding results for the case with the analysis prior. In the second set of experiments, ${\bf W}$ is an orthogonal Haar wavelet basis; the results are reported in Table~\ref{tab:deconvl1orthogonal} for the synthesis case, and in Table~\ref{tab:deconvl1orthogonal_analysis} for the analysis case. To visually illustrate the relative speed of the algorithms, Figure~\ref{fig:evolutioncriterionl1_redundant} plots the evolution of the constraint $\|{\bf B u}_k - {\bf y}\|$,
versus time, in experiments $1$,
for the synthesis prior case, with redundant wavelets.

\begin{table*}[t]
\centering \caption{Image deblurring using wavelets (redundant) - Computational load}
\label{tab:deconvl1redundant}
\begin{tabular}{|l|c|c|c| c|c|c| c|c|c|}
\hline
\footnotesize Expt. & \multicolumn{3}{|c|}{\footnotesize Avg. calls to ${\bf B},{\bf B}^H$ (min/max)} & \multicolumn{3}{|c|}{\footnotesize Iterations} & \multicolumn{3}{|c|}{\footnotesize CPU time (seconds)}\\
\hline
& \footnotesize SPGL1 & \footnotesize NESTA &\footnotesize  C-SALSA & \footnotesize SPGL1 & \footnotesize NESTA &\footnotesize  C-SALSA & \footnotesize SPGL1 & \footnotesize NESTA &\footnotesize  C-SALSA\\
\hline
\footnotesize 1 & \footnotesize 1029 (659/1290) & \footnotesize 3520 (3501/3541) & \footnotesize 398 (388/406)   & \footnotesize 340 & \footnotesize 880 & \footnotesize 134    & \footnotesize 441.16 & \footnotesize 590.79 & \footnotesize 100.72 \\
\footnotesize 2A & \footnotesize 511 (279/663) & \footnotesize 4897 (4777/4981) & \footnotesize 451 (442/460)    & \footnotesize 160 & \footnotesize 1224 & \footnotesize 136   & \footnotesize 202.67 & \footnotesize 798.81 & \footnotesize 98.85 \\
\footnotesize 2B & \footnotesize 377 (141/532) & \footnotesize 3397 (3345/3473) & \footnotesize 362 (355/370)    & \footnotesize 98 & \footnotesize 849 & \footnotesize 109     & \footnotesize 120.50 & \footnotesize 557.02 & \footnotesize 81.69 \\
\footnotesize 3A & \footnotesize 675 (378/772) & \footnotesize 2622 (2589/2661) & \footnotesize 172 (166/175)    & \footnotesize 235 & \footnotesize 656 & \footnotesize 58     & \footnotesize 266.41 & \footnotesize 423.41 & \footnotesize 42.56 \\
\footnotesize 3B & \footnotesize 404 (300/475) & \footnotesize 2446 (2401/2485) & \footnotesize 134 (130/136)    & \footnotesize 147 & \footnotesize 551 & \footnotesize 41     & \footnotesize 161.17 & \footnotesize 354.59 & \footnotesize 29.57 \\
\hline
\end{tabular}
\end{table*}

\begin{table*}
\centering \caption{Image Deconvolution using wavelets (redundant, analysis prior) - Computational load}
\label{tab:deconvl1redundant_analysis}
\begin{tabular}{|l|c|c| c|c|c|c|}
\hline
\footnotesize Expt. & \multicolumn{2}{|c|}{\footnotesize Avg. calls to ${\bf B},{\bf B}^H$ (min/max)} & \multicolumn{2}{|c|}{\footnotesize Iterations} & \multicolumn{2}{|c|}{\footnotesize CPU time (seconds)}\\
\hline
 & \footnotesize NESTA &\footnotesize  C-SALSA & \footnotesize NESTA &\footnotesize  C-SALSA & \footnotesize NESTA &\footnotesize  C-SALSA\\
\hline
\footnotesize 1 & \footnotesize 2881 (2861/2889) & \footnotesize 413 (404/419)  & \footnotesize 720 & \footnotesize 138     & \footnotesize 353.88 & \footnotesize 80.32 \\
\footnotesize 2A & \footnotesize 2451 (2377/2505) & \footnotesize 362 (344/371)     & \footnotesize 613 & \footnotesize 109     & \footnotesize 291.14 & \footnotesize 62.65 \\
\footnotesize 2B & \footnotesize 2139 (2065/2197) & \footnotesize 290 (278/299)     & \footnotesize 535 & \footnotesize 87  & \footnotesize 254.94 & \footnotesize 50.14 \\
\footnotesize 3A & \footnotesize 2203 (2181/2217) & \footnotesize 137 (134/143)     & \footnotesize 551 & \footnotesize 42  & \footnotesize 261.89 & \footnotesize 23.83 \\
\footnotesize 3B & \footnotesize 1967 (1949/1985) & \footnotesize 116 (113/119)     & \footnotesize 492 & \footnotesize 39  & \footnotesize 236.45 & \footnotesize 22.38 \\
\hline
\end{tabular}
\end{table*}

\begin{table*}
\centering \caption{Image deblurring using wavelets (orthogonal) - Computational load}
\label{tab:deconvl1orthogonal}
\begin{tabular}{|l|c|c|c| c|c|c| c|c|c|}
\hline
\footnotesize Expt. & \multicolumn{3}{|c|}{\footnotesize Avg. calls to ${\bf B},{\bf B}^H$ (min/max)} & \multicolumn{3}{|c|}{\footnotesize Iterations} & \multicolumn{3}{|c|}{\footnotesize CPU time (seconds)}\\
\hline
& \footnotesize SPGL1 & \footnotesize NESTA &\footnotesize  C-SALSA & \footnotesize SPGL1 & \footnotesize NESTA &\footnotesize  C-SALSA & \footnotesize SPGL1 & \footnotesize NESTA &\footnotesize  C-SALSA\\
\hline
\footnotesize 1 & \footnotesize 730 (382/922) & \footnotesize 13901 (13871/13931) & \footnotesize 494 (424/748)      & \footnotesize 298 & \footnotesize 3475 & \footnotesize 166   & \footnotesize 46.64 & \footnotesize 622.09 & \footnotesize 23.91 \\
\footnotesize 2A & \footnotesize 352 (191/480) & \footnotesize 1322 (1301/1329) & \footnotesize 205 (202/205)    & \footnotesize 128 & \footnotesize 331 & \footnotesize 69     & \footnotesize 19.21 & \footnotesize 58.29 & \footnotesize 10.07 \\
\footnotesize 2B & \footnotesize 207 (128/254) & \footnotesize 1218 (1193/1261) & \footnotesize 123 (115/133)    & \footnotesize 87 & \footnotesize 305 & \footnotesize 42  & \footnotesize 12.23 & \footnotesize 52.92 & \footnotesize 6.35 \\
\footnotesize 3A & \footnotesize 248 (161/320) & \footnotesize 1421 (1413/1433) & \footnotesize 118 (115/121)    & \footnotesize 104 & \footnotesize 355 & \footnotesize 40     & \footnotesize 14.98 & \footnotesize 58.693 & \footnotesize 5.57 \\
\footnotesize 3B& \footnotesize 170 (114/220) & \footnotesize 4408 (4345/4545) & \footnotesize 258 (94/328)      & \footnotesize 72 & \footnotesize 1102 & \footnotesize 87     & \footnotesize 9.51 & \footnotesize 181.83 & \footnotesize 11.93 \\
\hline
\end{tabular}
\end{table*}

\begin{table*}
\centering \caption{Image deblurring using wavelets (orthogonal, analysis prior) - Computational load}
\label{tab:deconvl1orthogonal_analysis}
\begin{tabular}{|l|c|c| c|c|c|c|}
\hline
\footnotesize Expt. & \multicolumn{2}{|c|}{\footnotesize Avg. calls to ${\bf B},{\bf B}^H$ (min/max)} & \multicolumn{2}{|c|}{\footnotesize Iterations} & \multicolumn{2}{|c|}{\footnotesize CPU time (seconds)}\\
\hline
 & \footnotesize NESTA &\footnotesize  C-SALSA & \footnotesize NESTA &\footnotesize  C-SALSA & \footnotesize NESTA &\footnotesize  C-SALSA\\
\hline
\footnotesize 1  & \footnotesize 8471 (8413/8553) & \footnotesize 387 (380/395)     & \footnotesize 2118 & \footnotesize 117    & \footnotesize 300.60 & \footnotesize 16.51 \\
\footnotesize 2A & \footnotesize 2463 (2445/2489) & \footnotesize 377 (371/383)     & \footnotesize 616 & \footnotesize 126     & \footnotesize 311.49 & \footnotesize 77.75 \\
\footnotesize 2B & \footnotesize 2159 (2097/2253) & \footnotesize 300 (290/317)     & \footnotesize 540 & \footnotesize 101     & \footnotesize 280.35 & \footnotesize 59.75 \\
\footnotesize 3A & \footnotesize 2203 (2165/2229) & \footnotesize 153 (149/155)     & \footnotesize 551 & \footnotesize 52  & \footnotesize 282.12 & \footnotesize 32.02 \\
\footnotesize 3B & \footnotesize 4710 (4577/4829) & \footnotesize 212 (104/374)     & \footnotesize 1178 & \footnotesize 59     & \footnotesize 167.73 & \footnotesize 7.89 \\
\hline
\end{tabular}
\end{table*}

\subsection{Image Deblurring with Total Variation}\label{sec:TVreg_exp}

The same five image deconvolution problems listed in Table~\ref{decon_problems} were also addressed using total variation (TV) regularization (more specifically, the isotropic discrete total variation, as defined in \cite{Chambolle}). The corresponding Moreau proximal mapping is computed using $5$ iterations of Chambolle's algorithm \cite{Chambolle}.

\begin{figure}
\centering
\subfigure[]{
\includegraphics[width=0.3\textwidth]{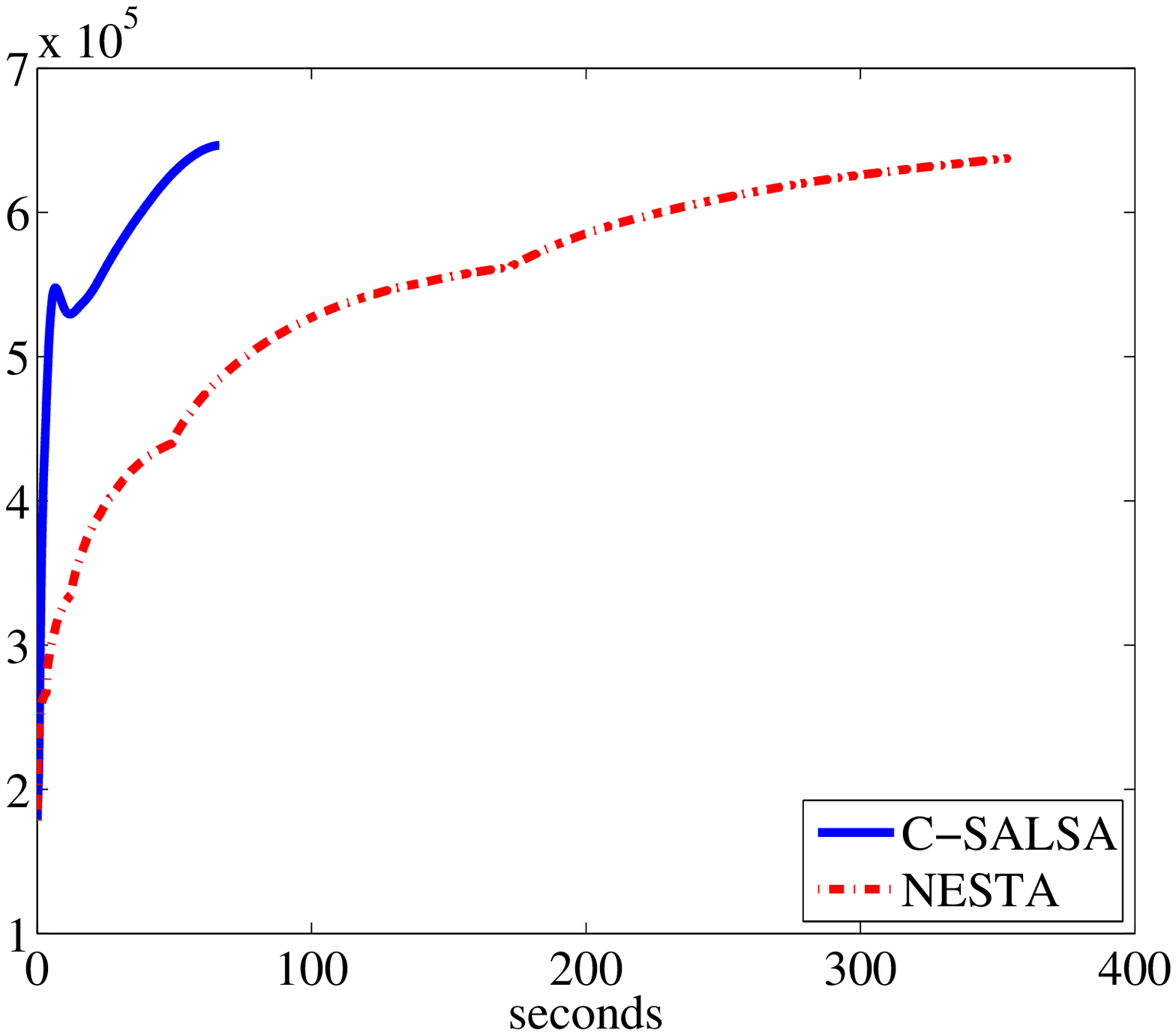}
}
\subfigure[]{
\includegraphics[width=0.3\textwidth]{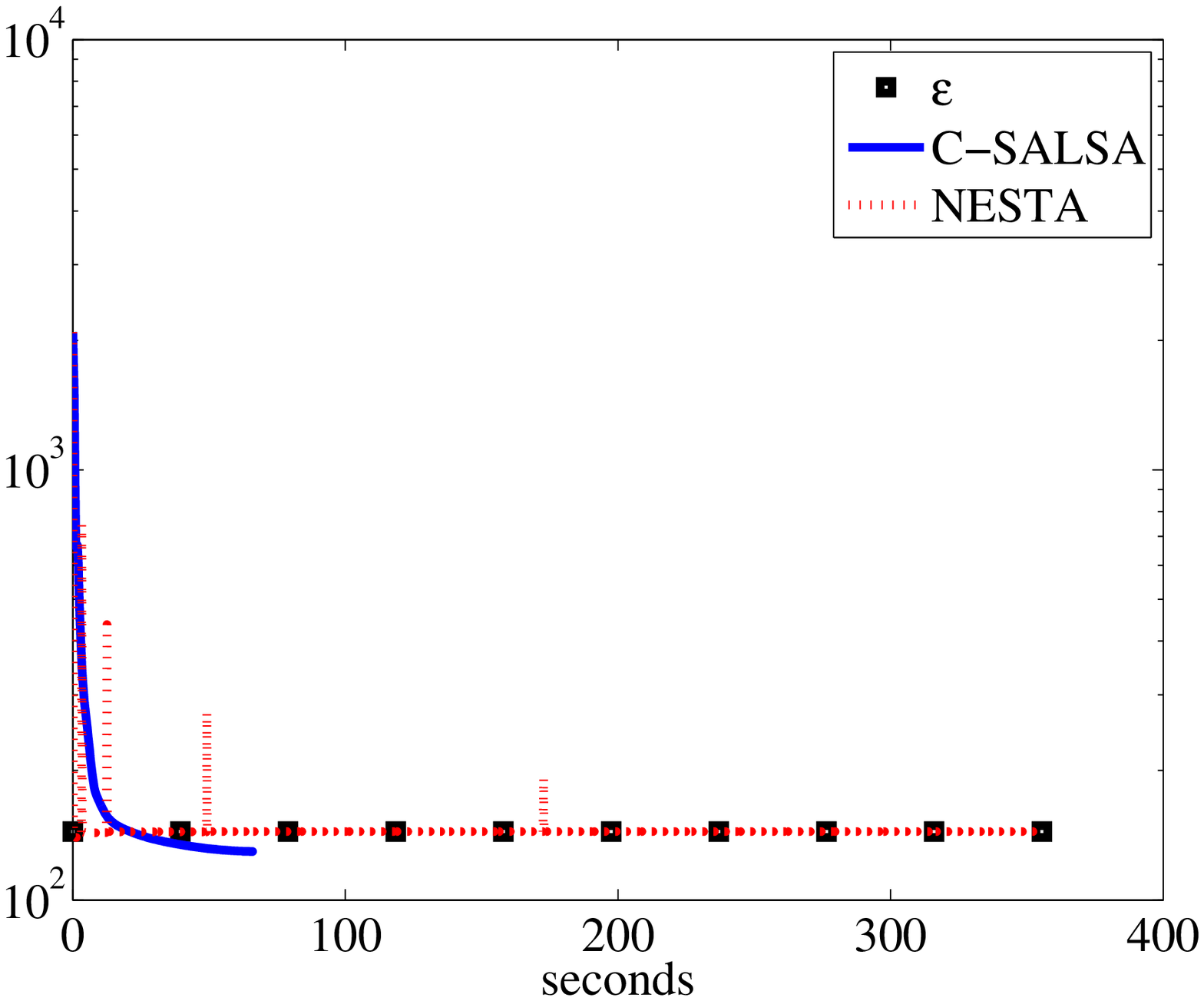}
}
\caption{Image deblurring (uniform blur) with TV regularization: (a) evolution of the objective function over time; (b) evolution of the constraint $\|{\bf B x}- {\bf y}\|$ over time.}
\label{fig:plots_deblurringTV}
\end{figure}

Table~\ref{tab:deconvTV} compares the performance of C-SALSA and NESTA, in terms of speed. The evolutions of the objective function and the constraint for experiment $1$ are plotted in Figure~\ref{fig:plots_deblurringTV}.

\begin{table*}
\centering \caption{Image deblurring using TV - Computational load}
\label{tab:deconvTV}
\begin{tabular}{|l|c|c| c|c|c|c|}
\hline
\footnotesize Expt. & \multicolumn{2}{|c|}{\footnotesize Avg. calls to ${\bf B},{\bf B}^H$ (min/max)} & \multicolumn{2}{|c|}{\footnotesize Iterations} & \multicolumn{2}{|c|}{\footnotesize CPU time (seconds)}\\
\hline
 & \footnotesize NESTA &\footnotesize  C-SALSA & \footnotesize NESTA &\footnotesize  C-SALSA & \footnotesize NESTA &\footnotesize  C-SALSA\\
\hline
\footnotesize 1 & \footnotesize 7783 (7767/7795) & \footnotesize 695 (680/710)  & \footnotesize 1945 & \footnotesize 232    & \footnotesize 311.98 & \footnotesize 62.56 \\
\footnotesize 2A & \footnotesize 7323 (7291/7351) & \footnotesize 559 (536/578)     & \footnotesize 1830 & \footnotesize 150    & \footnotesize 279.36 & \footnotesize 38.63 \\
\footnotesize 2B & \footnotesize 6828 (6775/6883) & \footnotesize 299 (269/329)     & \footnotesize 1707 & \footnotesize 100    & \footnotesize 265.35 & \footnotesize 25.47 \\
\footnotesize 3A & \footnotesize 6594 (6513/6661) & \footnotesize 176 (98/209)  & \footnotesize 1649 & \footnotesize 59     & \footnotesize 250.37 & \footnotesize 15.08 \\
\footnotesize 3B & \footnotesize 5514 (5417/5585) & \footnotesize 108 (104/110)     & \footnotesize 1379 & \footnotesize 37     & \footnotesize 210.94 & \footnotesize 9.23 \\
\hline
\end{tabular}
\end{table*}

We can conclude from Tables \ref{tab:deconvl1redundant}, \ref{tab:deconvl1redundant_analysis}, \ref{tab:deconvl1orthogonal}, \ref{tab:deconvl1orthogonal_analysis}, and \ref{tab:deconvTV} that, in image deconvolution problems, both with wavelet-based and TV-based
regularization, C-SALSA is almost always clearly faster than the fastest of the other competing algorithms.

\subsection{MRI Image Reconstruction}

\begin{figure}[h]
\centering
\subfigure[]{
\includegraphics[width=0.25\textwidth]{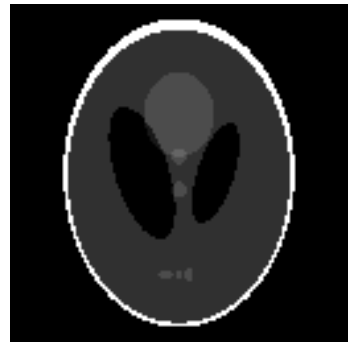}
\label{fig:phantom128}
}

\subfigure[]{
\includegraphics[width=0.25\textwidth]{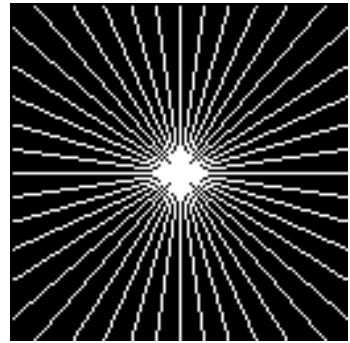}
\label{fig:MRImask22}
}

\subfigure[]{
\includegraphics[width=0.25\textwidth]{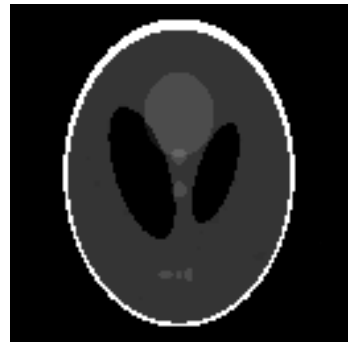}
\label{fig:estimateSALSA_mri}
}
\caption{MRI reconstruction: (a)$128\times 128$ Shepp Logan phantom; (b) Mask with 22 radial lines;  (c) image estimated using C-SALSA.}
\end{figure}

We now consider the problem of reconstructing the $128\times 128$ Shepp-Logan phantom (shown in Figure~\ref{fig:phantom128}) from a limited number of radial lines (22, in our experiments, as shown in Figure~\ref{fig:MRImask22}) of its 2D discrete Fourier transform.  The projections are also corrupted with circular complex Gaussian noise, with variance $\sigma^2\ =\ 0.5\times 10^{-6}$. We use TV regularization (as described in Subsection \ref{sec:TVreg_exp}), with the corresponding Moreau proximal mapping implemented by  $10$ iterations of Chambolle's algorithm \cite{Chambolle}.

\begin{table}
\centering \caption{MRI reconstruction - Comparison}
\label{tab:MRIresults}
\begin{tabular}{|l|c| c|c|c|}
\hline
\footnotesize Algorithm & \footnotesize Calls to ${\bf B},{\bf B}^H$ & \footnotesize Iterations & \footnotesize time (seconds) & \footnotesize MSE\\
\hline
\footnotesize NESTA & \footnotesize 1228 (1161/1261) & \footnotesize 307 & \footnotesize 15.50 & \footnotesize 9.335e-6 \\
\hline
\footnotesize C-SALSA & \footnotesize 366 (365/368) & \footnotesize 122 & \footnotesize 12.89 & \footnotesize 2.440e-6 \\
\hline
\end{tabular}
\end{table}

\begin{figure}[th]
\centering
\subfigure[]{
\includegraphics[width=0.3\textwidth]{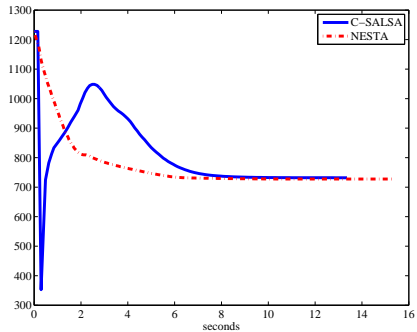}
}
\subfigure[]{
\includegraphics[width=0.3\textwidth]{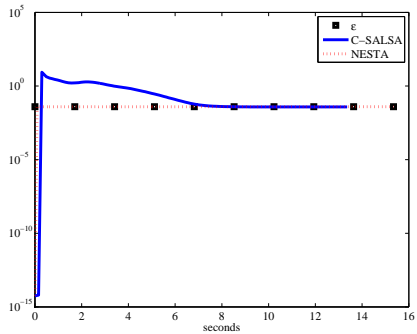}
}

\caption{MRI reconstruction with TV regularization: (a) evolution of the objective function over time; (b) evolution of the constraint $\|{\bf B x}- {\bf y}\|$ over time.}
\label{fig:plots_mri}
\end{figure}

Table~\ref{tab:MRIresults} shows the number of calls, number of iterations, and CPU times, while Figure~\ref{fig:plots_mri} plots
the evolution of the objective function and constraint over time.
Figure~\ref{fig:estimateSALSA_mri} shows the estimate obtained using C-SALSA (the estimate NESTA  is, naturally,
visually indistinguishable). Again, we may conclude that C-SALSA is faster than NESTA, while achieving comparable values of mean squared error of
the reconstructed image.

\subsubsection*{High Dynamic Range TV Reconstruction}
A related example that we will consider here is the reconstruction
of images composed of random squares, from their partial Fourier
measurements, with TV regularization (see \cite{NESTA}, section
$6.4$). The dynamic range of the signals (the amplitude of the
squares) varies from 20 dB to 80 dB. The size of each image is $128
\times 128$, the number of radial lines in the DFT measurement mask
is $27$ (corresponding to $m / n \approx 0.2$), and the Gaussian
noise standard deviation is $\sigma = 0.1$.

\begin{figure}[h!]
\centering
\label{fig:squares_hdr}
\includegraphics[width=0.35\textwidth]{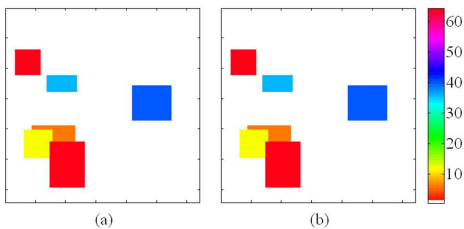}
\caption{TV image reconstruction: (a) Original image with dynamic range $=40$ dB; (b) Estimate using C-SALSA.}
\end{figure}

\begin{figure}[th]
\centering
\subfigure[]{
\includegraphics[width=0.3\textwidth]{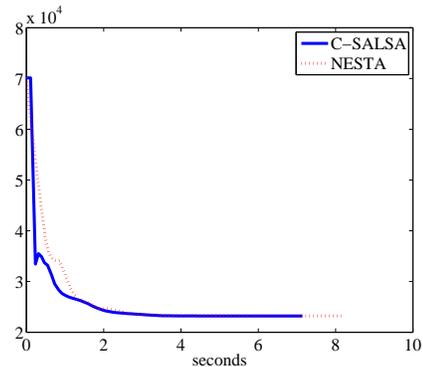}
}
\subfigure[]{
\includegraphics[width=0.3\textwidth]{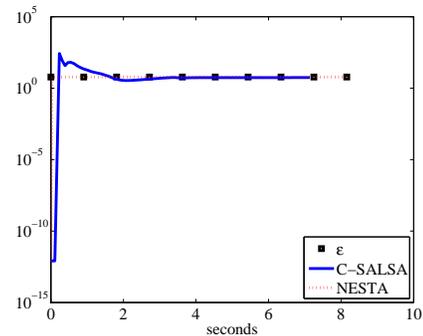}
}
\caption{TV image reconstruction (dynamic range $= 40$ dB): (a) evolution of the objective function over time; (b) evolution of the constraint $\|{\bf B x}- {\bf y}\|$ over time.}
\label{fig:plots_squares}
\end{figure}

Figure~\ref{fig:squares_hdr} shows the original image with a dynamic range of $40$ dB and the estimate obtained using C-SALSA. Figure~\ref{fig:plots_squares} shows the evolution over time of the objective and the error constraint for C-SALSA and NESTA, while Table~\ref{tab:squares_hdr} compares the two algorithms with respect to the number of calls to ${\bf A}$ and ${\bf A}^H$, number of iterations, CPU time, and MSE obtained, over 10 random trials. It is clear from Table~\ref{tab:squares_hdr} that C-SALSA uses considerably fewer calls to the operators ${\bf A}$ and ${\bf A}^H$ than NESTA.

\begin{table*}
\centering \caption{Image reconstruction (high dynamic range) using TV - Computational load}
\label{tab:squares_hdr}
\begin{tabular}{|l|c|c| c|c|c|c|c|c|}
\hline
\footnotesize Dyn. range & \multicolumn{2}{|c|}{\footnotesize Avg. calls to ${\bf B},{\bf B}^H$ (min/max)} & \multicolumn{2}{|c|}{\footnotesize Iterations} & \multicolumn{2}{|c|}{\footnotesize CPU time (seconds)} & \multicolumn{2}{|c|}{\footnotesize MSE}\\
\hline
\footnotesize (dB) & \footnotesize NESTA &\footnotesize  C-SALSA & \footnotesize NESTA &\footnotesize  C-SALSA & \footnotesize NESTA &\footnotesize  C-SALSA & \footnotesize NESTA &\footnotesize  C-SALSA\\
\hline
\footnotesize 20 & \footnotesize 1213 (1169/1273) & \footnotesize 226 (224/227) & \footnotesize 303 & \footnotesize 76 & \footnotesize 8.99& \footnotesize 7.24 & \footnotesize 0.00241743 & \footnotesize 0.000543426 \\
\footnotesize 40 & \footnotesize 991 (961/1017) & \footnotesize 227 (224/227) & \footnotesize 248 & \footnotesize 76 & \footnotesize 7.34 & \footnotesize 7.002 & \footnotesize 0.00432206 & \footnotesize 0.000651107 \\
\footnotesize 60 & \footnotesize 731 (721/737) & \footnotesize 282 (281/284) & \footnotesize 183 & \footnotesize 95 & \footnotesize 4.92 & \footnotesize 8.35 & \footnotesize 0.005294 & \footnotesize 0.00072848 \\
\footnotesize 80 & \footnotesize 617 (613/617) & \footnotesize 353 (350/353) & \footnotesize 154 & \footnotesize 118 & \footnotesize 4.16 & \footnotesize 10.72 & \footnotesize 0.00702862 & \footnotesize 0.000664638 \\
\hline
\end{tabular}
\end{table*}

\subsection{Image Inpainting}
\begin{figure}[h!]
\centering
\subfigure[]{
\includegraphics[width=0.25\textwidth]{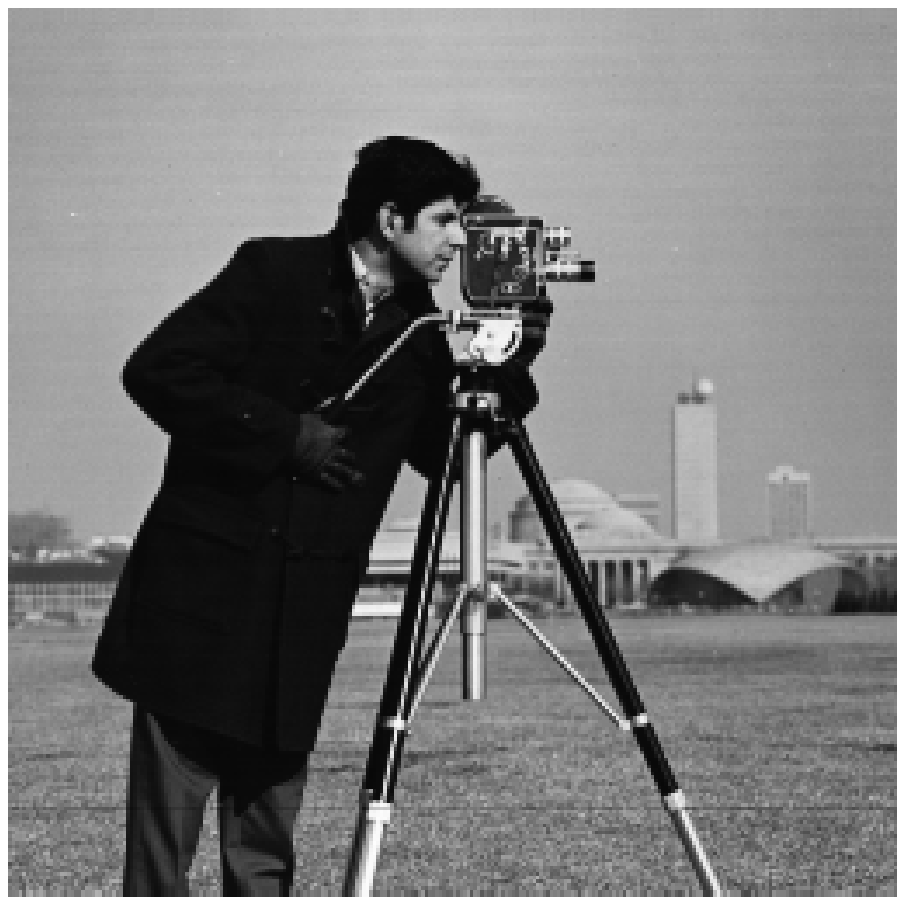}
}
\subfigure[]{
\includegraphics[width=0.25\textwidth]{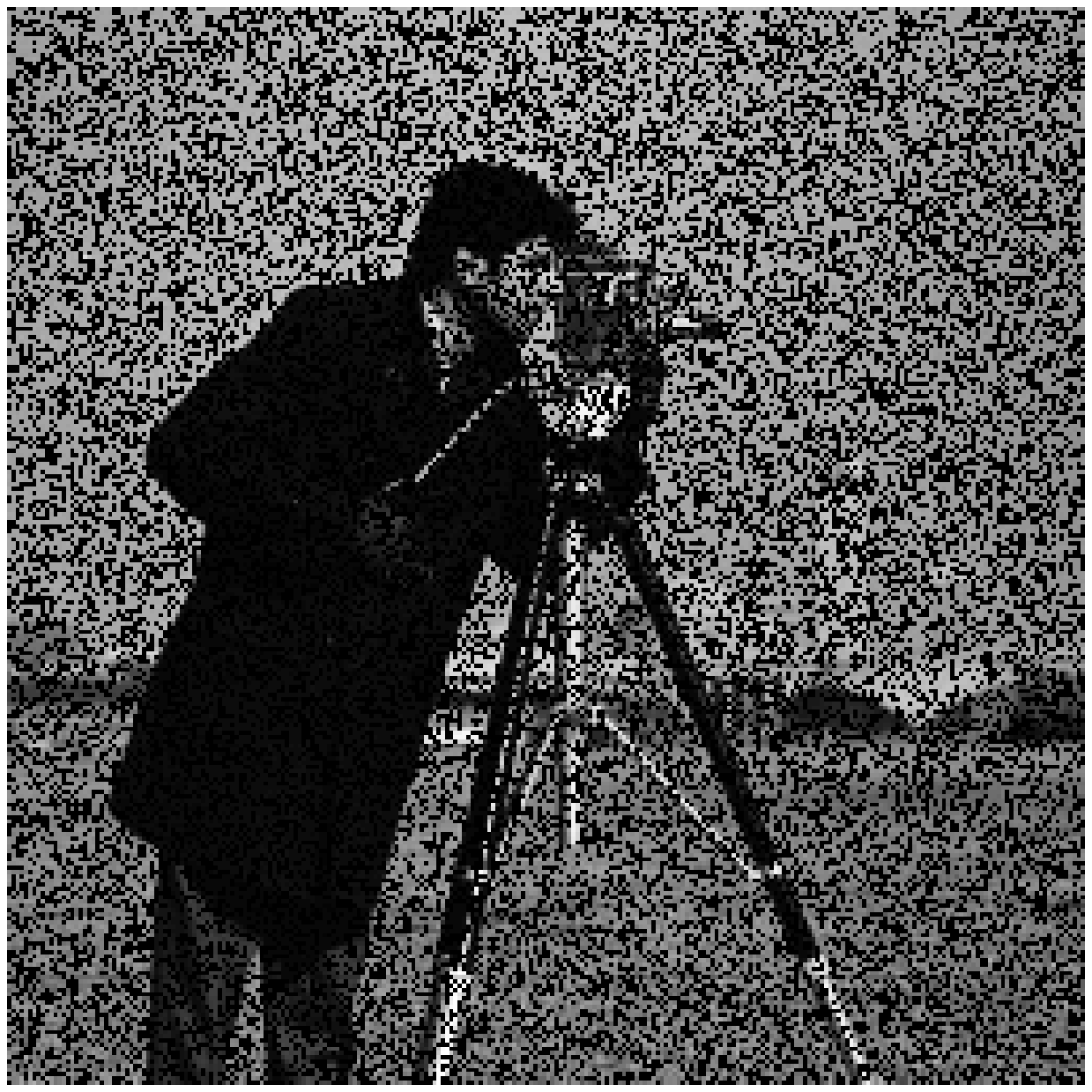}
}
\subfigure[]{
\includegraphics[width=0.25\textwidth]{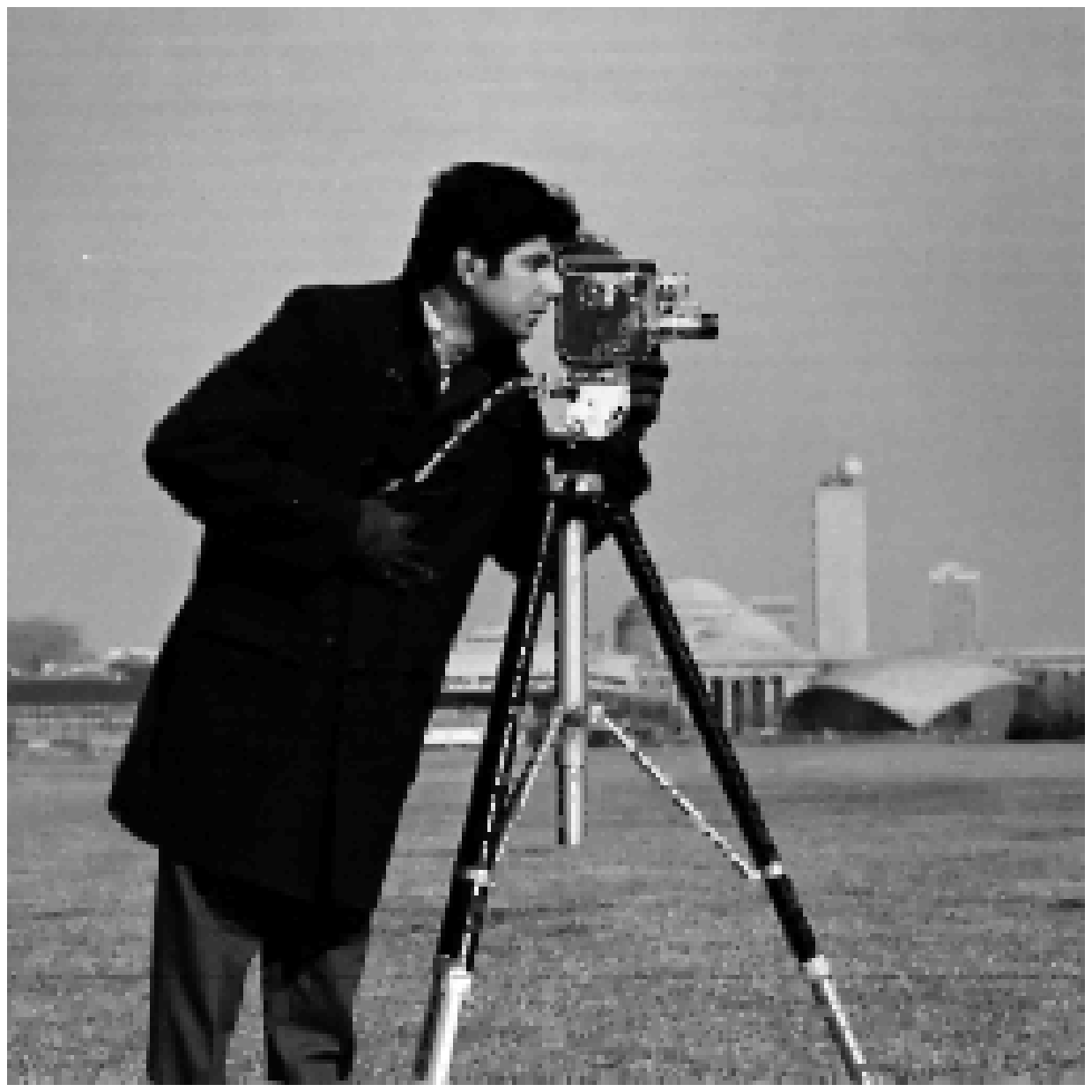}
}
\caption{Image inpainting with TV regularization: (a) Original cameraman image; (b) Image with $40\%$ pixels missing; (c) Estimated using C-SALSA.}
\label{fig:missingpixels}
\end{figure}

Finally, we consider an image inpainting problem, as explained in Section~\ref{sec:computingxk}.
The original image is again the Cameraman, and the observation consists in losing $40\%$ of its pixels,
as shown in Figure~\ref{fig:missingpixels}. The observations are also corrupted with Gaussian noise
(with an SNR of $40$ dB). The regularizer is again TV, implemented by $10$ iterations of Chambolle's algorithm.

\begin{figure}[th]
\centering
\subfigure[]{
\includegraphics[width=0.3\textwidth]{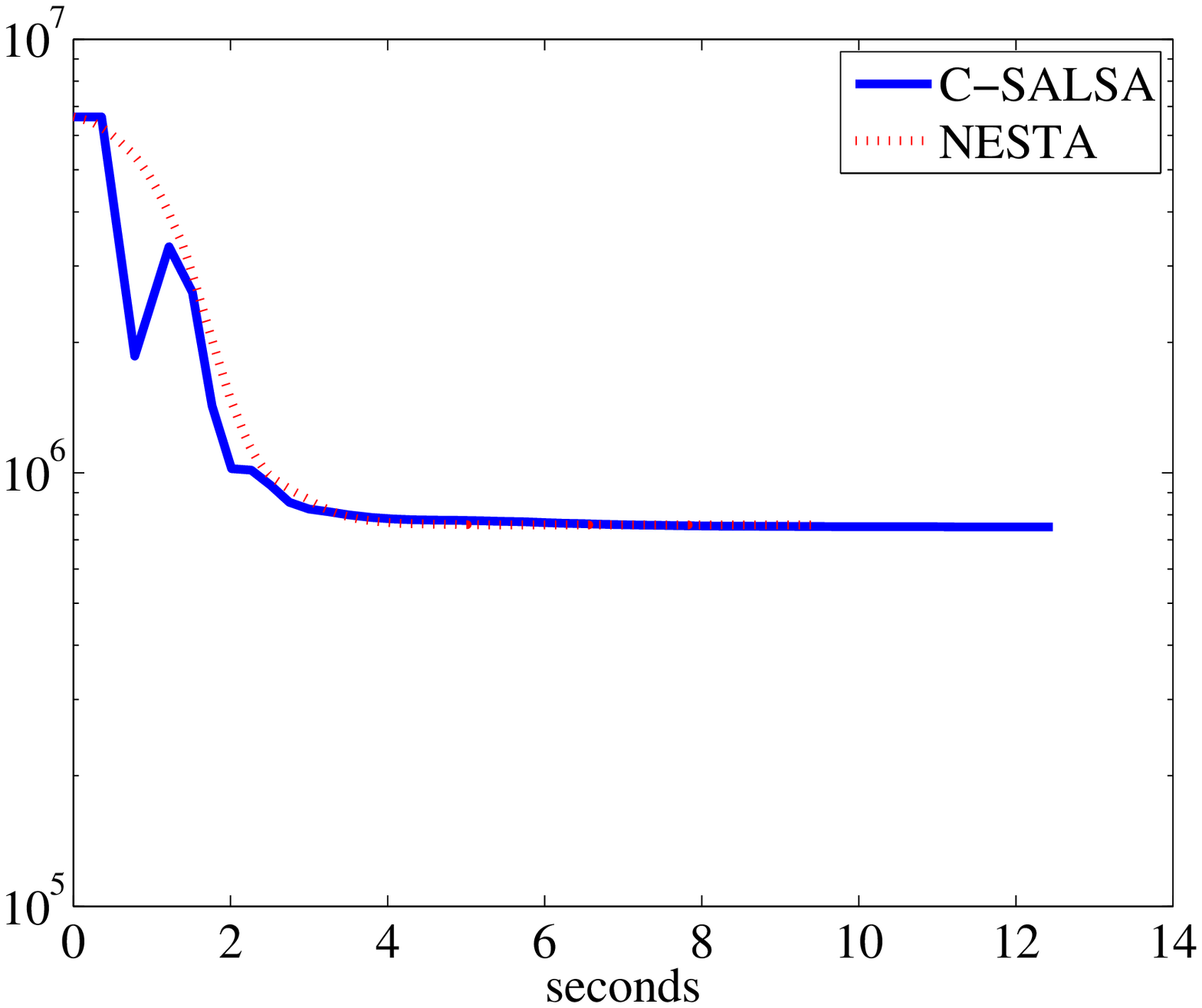}
}
\subfigure[]{
\includegraphics[width=0.3\textwidth]{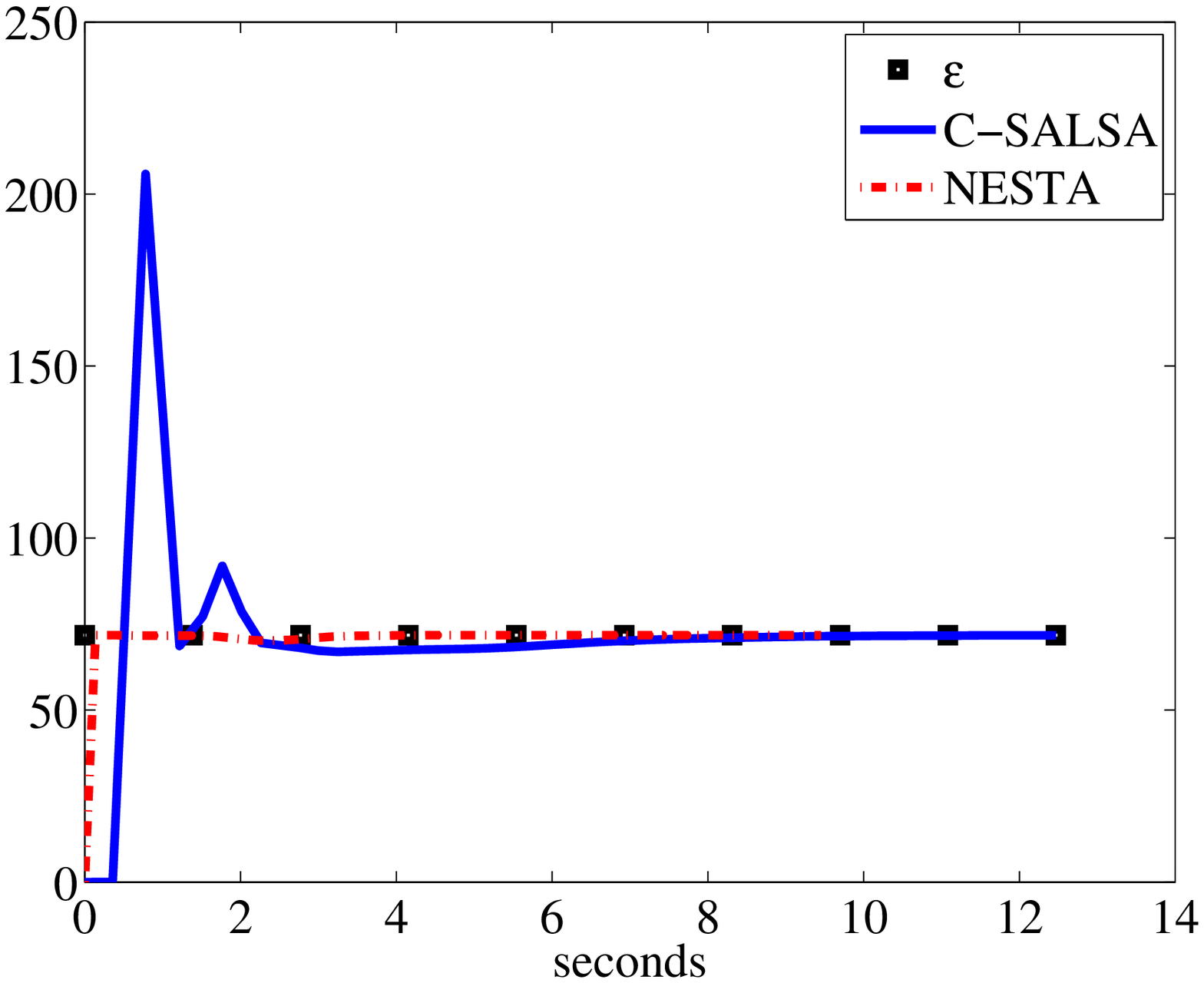}
}
\caption{Image inpainting with TV regularization: (a) evolution of the objective function over time; (b) evolution of the constraint $\|{\bf B x}- {\bf y}\|$ over time.}
\label{fig:objective_missingpixels}
\end{figure}

The image estimate obtained by C-SALSA is shown in Figure~\ref{fig:missingpixels}, with the original
also shown for comparison. The estimate obtained using NESTA was visually very similar. Table~\ref{tab:missingdata_comparison} compares the performance of the two algorithms, and Figure~\ref{fig:objective_missingpixels} shows the evolution of the objective function for each of them.

\begin{table}[hbt]
\centering \caption{Image inpainting: Comparison.}
\label{tab:missingdata_comparison}
\begin{tabular}{|c|l|l|l|l|l|}
  \hline
  & \footnotesize Calls to ${\bf B},{\bf B}^H$ & \footnotesize Iterations & \footnotesize time (seconds) & \footnotesize MSE \\
\hline
\footnotesize NESTA & \footnotesize 403 (401/405) & \footnotesize 101 & \footnotesize 10.29 & \footnotesize 81.316 \\
\hline
\footnotesize C-SALSA & \footnotesize 143 (143/143) & \footnotesize 47 & \footnotesize 12.97 & \footnotesize 75.003 \\
\hline

\end{tabular}
\end{table}

\Section{Conclusions}
\label{sec:conclusions}

We have presented a new algorithm for solving the constrained
optimization formulation of regularized image reconstruction/restoration.
The approach, which can be used with any type of convex regularizers
(wavelet-based, total variation), is based on a VS technique
which yields an equivalent constrained problem. This constrained
problem is then addressed using an augmented Lagrangian method,
more specifically, the alternating direction method of multipliers (ADMM).
Our algorithm works for any convex regularizer for which the Moreau proximal mapping can be efficiently computed, and is therefore more general purpose than some of the available state-of-the-art methods which are available only for either $\ell_1$- and/or TV regularization. Experiments on a set of standard image recovery problems (deconvolution, MRI reconstruction, inpainting) have shown that the proposed algorithm (termed C-SALSA, for {\it constrained split augmented Lagrangian shrinkage algorithm}) is usually faster than previous state-of-the-art methods. Automating the choice of the value of the parameter $\mu$ remains an open question.

\end{document}